\documentclass[10pt]{article}
\usepackage[english]{babel}
 \pdfoutput=1 
\usepackage{amsmath,amssymb,mathrsfs}
\usepackage{pgf}

\usepackage[latin1]{inputenc}
\usepackage{color}
\usepackage{epsfig}
\usepackage{psfrag}
\usepackage{eucal}
\usepackage{paralist}

\usepackage{latexsym,amssymb}
\usepackage{hyperref}
\usepackage{amsmath,amsthm,amsxtra}
\usepackage{amsfonts}
\usepackage{epsfig,inputenc}
\usepackage{graphpap,latexsym,epsf}
\usepackage{color}
\usepackage{amssymb, paralist,color,enumerate}
\usepackage{graphicx}
\usepackage{verbatim}

\usepackage{bbm} 

\usepackage{pgf}

\definecolor{dred}{rgb}{.8,0,0}
\definecolor{ddmagenta}{rgb}{0.7,0,0.9}
\definecolor{ddcyan}{rgb}{0,0.2,1.0}

\newtheorem{theorem}{Theorem}[section]

\newtheorem{lemma}[theorem]{Lemma}

\newtheorem{proposition}[theorem]{Proposition}
\newtheorem{examples}[theorem]{Examples}
\newtheorem{example}[theorem]{Example}

\newtheorem{corollary}[theorem]{Corollary}

\newtheorem{remark}[theorem]{Remark}
\newtheorem{definition}[theorem]{Definition}

\numberwithin{equation}{section}
 
\newcommand{\QED}{\mbox{}\hfill\rule{5pt}{5pt}\medskip\par}
\newcommand{\dd}{\, \mathrm{d}}

\numberwithin{equation}{section}
\numberwithin{figure}{section}

\newcommand{\Nz}{{\mathbb{N}}}

\newcommand\JUMP[1]{\mathchoice
                  {\big[\hspace*{-.3em}\big[#1\big]\hspace*{-.3em}\big]}
                   {[\hspace*{-.15em}[#1]\hspace*{-.15em}]}
                   {[\![#1]\!]}
                   {[\![#1]\!]}}

\newcommand{\supp}{\mathop{\mathrm{supp}}}

\def\calD{{\mathcal D}} \def\calE{{\mathcal E}} \def\calF{{\mathcal F}}
  
\def\calJ{{\mathcal J}}  \def\calL{{\mathcal L}}
  
  \def\calR{{\mathcal R}}

\newcommand{\Spz}{\mathbf{Z}}
\newcommand{\Spx}{\mathbf{X}}

\newcommand\GC{\Gamma_{\mbox{\tiny\rm C}}}

\newcommand{\N}{\mathbb{N}}
\newcommand{\R}{\mathbb{R}}
\newcommand{\eps}{\varepsilon}

\newcommand{\aein}{\text{a.e. in}}
\newcommand{\foraa}{\text{for a.a.}}

\newcommand{\ddx}{\,\mathrm{d}}

\newcommand{\piecewiseConstant}[2]{\overline{#1}_{\kern-1pt#2}}
\newcommand{\pwc}{\piecewiseConstant}
\newcommand{\underpiecewiseConstant}[2]{\underline{#1}_{\kern-1pt#2}}
\newcommand{\upwc}{\underpiecewiseConstant}

\newcommand{\piecewiseLinear}[2]{#1_{\kern-1pt#2}}
\newcommand{\pwl}{\piecewiseLinear}

\newcommand{\weaksto}{\overset{*}{\rightharpoonup}}

\newcommand{\BV}{\mathrm{BV}}

\newcommand{\SBV}{\mathrm{SBV}}

\newenvironment{rcomm}{ \small  \color{ddcyan} \bf R: }{\normalsize \color{black}}
\newenvironment{rnew}{\color{dred}}{\normalsize \color{black}}

\newcommand{\berin}{\begin{rnew}}
\newcommand{\erin}{\end{rnew}}
\newcommand{\beroc}{\begin{rcomm}}
\newcommand{\eroc}{\end{rcomm}}

\newcommand{\meanbar}[1]{%
  \setbox0 = \hbox{$#1 \int$}
  \hbox to 0pt{%
    \thinspace
    \hskip 0.22\wd0
    \raise 0.5\ht0
    \hbox{%
      \lower 0.5\dp0
      \hbox{\rule{0.6\wd0}{\linethic\tauess}}
    }%
    \hss
  }%
}

\newcommand{\sft}{\mathsf{t}}
\newcommand{\GMMt}[3]{\mathrm{GMM}([#1,#2];\Spz,\calE_{#3},\calD;Z_0)}
\newcommand{\GMM}[1]{\mathrm{GMM}(\Spz,\calE_{#1},\calD;Z_0)}
\newcommand{\SMM}[1]{\mathrm{SMM}(\Spz,\calE_{#1},\calD;Z_0)}
\newcommand{\SMMk}{\mathrm{SMM}(\Spz,\calE_k,\calD;Z_0^{k})}
\newcommand{\SMMt}[3]{\mathrm{SMM}([#1,#2];\Spz,\calE_{#3},\calD;Z_0)}
\newcommand{\conv}{\overset{*}{\rightharpoonup}}

\newcommand{\down}{\downarrow}

\definecolor{ddmagenta}{rgb}{0.7,0,0.9}
\definecolor{Turk}{rgb}{0,0.7,0.4}

\newcommand{\NEW}{\color{black}}

\definecolor{dorogreen}{rgb}{0,0.8,0}
\definecolor{violet}{rgb}{0.4,0,0.9}

\definecolor{ddcyan}{rgb}{0,0.2,1.0}

\newcommand{\dissparam}{ a}

\begin{document}

\title{Rate-independent evolution of sets}

\author{Riccarda Rossi
\thanks{
DIMI, University of Brescia, 
Via Branze 38, I--25133 Brescia, Italy.
Email: {\ttfamily riccarda.rossi@unibs.it}.} 
\and
Ulisse Stefanelli\thanks{University of Vienna,
Oskar-Morgenstern-Platz 1, 1090 Wien, Austria.
Email: {\ttfamily  ulisse.stefanelli@univie.ac.at}
} 
\and
Marita
Thomas\thanks{Weierstrass Institute for Applied Analysis and
Stochastics, Mohrenstr.~39, 10117 Berlin, Germany.
Email: {\ttfamily marita.thomas@wias-berlin.de}}}

\date{}

\maketitle

 \centerline{\sl  Dedicated to Alexander Mielke on the occasion of his 60th birthday}
 
 \medskip





\begin{abstract}
 The goal of this work is to 
analyze a model for the rate-independent evolution of sets with finite perimeter. 
The evolution of the admissible sets is driven by that of  a given time-dependent set,  
which 
has to include the admissible sets and hence is to be understood as an 
external loading. The process is driven by the competition
between  perimeter minimization and minimization of volume
changes. 
\par
In  the mathematical modeling of  this process, we   distinguish the \emph{adhesive} case, in which the constraint 
that the (complement of) the  `external load' contains the evolving sets is penalized by a  term contributing to the driving energy functional, from the \emph{brittle} case,  enforcing this  constraint. 
The existence of \emph{Energetic solutions} for the adhesive system is
proved by passing to the limit in the associated time-incremental
minimization scheme. In the brittle case, this time-discretization
procedure gives rise to evolving sets satisfying the stability
condition, but it remains an open problem to additionally deduce
energy-dissipation balance 
in the time-continuous limit.  This  can be obtained under some
suitable quantification of data.  
 The properties of the brittle  evolution law are illustrated by numerical examples in two space dimensions. 
\end{abstract}

\paragraph{MSC 2010.} 
35A15, 
35R37, 
49Q10, 
74R10. 
\paragraph{Keywords and phrases. }Unidirectional evolution of sets by competition of perimeter and volume, minimizers of perimeter perturbed by  a nonsmooth functional, Minimizing Movements, stability, Energetic solutions.

%
%
\section{Introduction} 
The aim of this work is to introduce and analyze a notion of \emph{rate-independent} evolution for a 
set-valued function $Z: [0,T] \rightrightarrows \Omega$ (with $\Omega$ a bounded domain in $\R^d$),
whose evolution is triggered by that  of
another, given set-valued function  $F: [0,T]\rightrightarrows \Omega$, in the position of an external force,
through the constraint
\begin{equation}
\label{brittle-constr-intro}
Z(t) \cap F(t)  =  \emptyset \quad \text{for every } t \in [0,T].
\end{equation}
The evolution of  $Z$  is additionally ruled by the competition between  
the minimization of the perimeter 
and that of volume changes.  
\subsection*{Related models: brittle delamination and adhesive contact}
Our study is inspired and motivated by the modeling of delamination between two (elastic) 
bodies $O_+$ and  $O_- \subset \R^m$,  bonded along a prescribed contact surface  $\Gamma=\overline{O}_{+}\cap\overline{O}_{-}$  over
time interval $ (0,T)$.  Following the approach by \textsc{M.\
  Fr\'emond} \cite{Frem88CA,Fre02},   this process can be described in
terms of the temporal evolution of a phase-field type parameter, the
delamination variable $z: (0,T) \times \Gamma \to [0,1]$, which
represents   the fraction of fully effective molecular links in the bonding. Therefore, $z(t,x) =1$ ($z(t,x)=0$, respectively) means that the bonding is fully intact (completely broken)  at a given time instant $t\in [0,T]$ and in a given material point $x\in \Gamma$.  In models for \emph{brittle delamination}, the evolution of $z$ is coupled to that of the (small-strain) displacement variable $u: (0,T) \times O \to \R^m$ (with  $O: = O_+ \cup O_- $) through the so-called
\begin{equation}
\label{brittle-constraint-delam}
\text{brittle constraint} \qquad z \JUMP{u} =0 \quad \text{a.e.\ in } (0,T) \times \Gamma.
\end{equation}
In \eqref{brittle-constraint-delam}, $\JUMP{u} = u^{+}|_{\Gamma} -  u^{-}|_{\Gamma} $ ($u^+,\, u^-$ denoting the restrictions of $u$ to $O_+,\, O_-$, respectively) is the \emph{jump of $u$} across the interface $\Gamma$. Therefore, \eqref{brittle-constraint-delam} ensures the continuity of
the displacements, i.e.\ $\JUMP{u}=0$, in the (closure of the) set of points where (a portion
of) the bonding is still active, i.e.\ $z>0$.  In fact, \eqref{brittle-constraint-delam}   allows for displacement jumps only at
points where the bonding is completely broken, namely where $z=0$.
The set $\Gamma \setminus \supp(z)$,  where the displacements may  jump, can be thus understood as a \emph{crack set}; indeed, brittle delamination can be interpreted as a model for brittle fracture, along a \emph{prescribed} $(m{-}1)$-dimensional surface. 
\par
To the best of our knowledge, the analysis of the above described  brittle delamination model
 has been carried out only in the case this process is treated 
as  \emph{rate-independent}. Namely, the evolution of the internal
variable, responsible for the dissipation of energy, is governed by a
dissipation potential $\calR$ which is positively homogeneous of degree $1$, i.e.\ fulfilling $\calR(\lambda \dot z) = \lambda \calR(\dot z)$ for all $\lambda\geq 0$. In that case, the process can be mathematically modeled by means  of the general concept  of \emph{Energetic solution} to a rate-independent system, pioneered in \cite{MieThe99MMRI, MieThe04RIHM} (cf.\ also the parallel notion of \emph{quasistatic evolution} for brittle fracture, \cite{Francfort-Marigo98,DM-Toa2002}). The existence of Energetic solutions to the brittle delamination system was  proven in 
\cite{RoScZa09QDP} by passing to the limit, as the penalization parameter $k$  blows up, in the Energetic formulation for an approximate system. Therein,
 the brittle constraint \eqref{brittle-constraint-delam} was indeed penalized by the  
\begin{equation}
\label{adhesive-intro}
\text{adhesive contact term} \qquad \int_{\Gamma} k z|\JUMP{u}|^2 \dd \mathcal{H}^{m-1}(x), \qquad k \in (0,\infty)
\end{equation}
contributing to the  driving energy. For  the  corresponding rate-independent system,  referred to as   \emph{adhesive contact} system, the existence of Energetic solution dates back to \cite{KoMiRo006RIAD}.
The 
 passage from adhesive to brittle was also  studied 
in \cite{RosTho12ABDM} by  extending the existence of
Energetic solutions to the case in which the rate-independent
evolution of the delamination parameter $z$ is coupled to the
\emph{rate-dependent} evolution of the displacement, ruled by  a (no
longer quasi-static) momentum balance  with viscosity. In order to
overcome the analytical difficulties attached to the coupling of
rate-independent/rate-dependent behavior, in \cite{RosTho12ABDM} a
gradient regularization is advanced for the delamination variable 
of $\BV$-type. In fact,
(1)
 the constraint $z\in \{0,1\}$ was added to the model, making it closer to  Griffith-type model for crack evolution; (2)  the term $|\mathrm{D}z|(\Gamma)$  (namely, the total variation of the measure $\mathrm{D}z$) was added to the driving energy functional. 
Since $z$ is the characteristic function of the set $Z:= \{ x\in \Gamma\, : \ z(x) >0\}$, the gradient regularization 
$
|\mathrm{D}z|(\Gamma)$   coincides with the perimeter $ P(Z,\Gamma)  $ of  $ Z $  in $ \Gamma. $
\par
The rate-independent evolution of sets $Z: [0,T]\rightrightarrows \Omega$ studied in this paper can be understood as 
 an abstraction of the delamination process addressed in \cite{RoScZa09QDP, RosTho12ABDM}. Indeed,  suppose that  $\Omega = \Gamma$, and that  the temporal evolution of the set 
\begin{equation}
\label{F-delamination}
F(t) := \{ x\in \GC\, : \ \JUMP{u(t,x)}\neq 0 \}
\end{equation}
is \emph{given}. The sets $Z(t)$ correspond to the supports of the delamination variable and, when the set-valued mapping $F:[0,T]\rightrightarrows \Omega$ is given by \eqref{F-delamination}, \eqref{brittle-constr-intro} renders
 the brittle constraint \eqref{brittle-constraint-delam}. Therefore, hereafter we shall refer to 
 \eqref{brittle-constr-intro}  as  a (generalized) \emph{brittle constraint}, too.
 Like in \cite{RosTho12ABDM}, the energy functional driving the evolution of the sets $Z$ will feature their perimeter $P(Z,\Omega)$ in $\Omega$. 
  Again in analogy with delamination processes, in addition to the brittle case, in which \eqref{brittle-constr-intro} is enforced, we will also address the \emph{adhesive} case, in which  \eqref{brittle-constr-intro} is suitably penalized. 
\par
Our aim is threefold:
\begin{compactenum}
\item prove the existence of Energetic solutions for both the \emph{adhesive} and the \emph{brittle} evolution of sets;
\item investigate to what extent the fine geometric  properties proved in  \cite{RosTho12ABDM} for (the supports of)  Energetic solutions  to the adhesive and brittle delamination systems, carry over to this generalized setting;
\item
gain insight into the connection between this rate-independent evolution of sets, and the well-known mean curvature flow. 
\end{compactenum}
\subsection*{Solution concepts and our existence results}
Throughout the paper,  we will work under the (additional, w.r.t.\ the delamination case \eqref{F-delamination}) assumption that the `external load' $F: [0,T]\rightrightarrows \Omega$ is monotonically increasing, namely
\begin{equation}
\label{F-increasing}
F(s)\subset F(t) \qquad \text{ for all $0\leq s \leq t \leq T$.}
\end{equation}
Hence, in view of the constraint $Z(t) \subset F^c(t)$ from \eqref{brittle-constr-intro}, it will be natural to enforce for the set function $Z: [0,T] \rightrightarrows \Omega $ the opposite monotonicity property, 
viz.\ $Z(t) \subset Z(s)$  for all $0\leq s \leq t \leq T$.
\par
The
 adhesive and the brittle processes will be mathematically modeled by 
 a triple $(\Spz,\calE_k,\calD)_k$, 
 $k \in \N \cup \{\infty\}$,  with 
 \begin{compactitem}
 \item[-]
 $\Spz = \{ Z \subset \Omega  \, : \ \calL^d(Z)<\infty\}$ the state space;
 \item[-] the driving energy functional 
 \[
 \calE_\infty(t,Z):=P(Z,\Omega)+\calJ_\infty(t,Z) \qquad \text{with } 
\calJ_\infty(t,Z):=
\begin{cases}
0 &\text{ if } Z \cap F(t) =\emptyset,
\\
\infty & \text{ otherwise},
\end{cases}
 \]
 for the brittle process, while for the adhesive  process we will pose   
 \[
 \calE_k(t,Z):=P(Z,\Omega)+\calJ_k(t,Z)\,, \qquad \text{where the term } 
\calJ_k(t,Z):= \int_\Omega k  f(t) \chi_Z \dd x, \ k \in (0,\infty),
 \] 
 penalizes the brittle constraint \eqref{brittle-constr-intro}, in that the support of  the  function $f(t)$ coincides with $F(t)$;
 \item[-] the dissipation (quasi-)distance 
 \[
 \calD:\Spz\times\Spz\to[0,\infty],\quad 
\calD(Z_1,Z_2):=\left\{
\begin{array}{ll}
\dissparam \calL^d(Z_1\setminus Z_2) &\text{if }Z_2\subset Z_1,\\
\infty&\text{otherwise}
\end{array}
\right.
 \]
 for some $\dissparam>0$, 
 enforcing that $Z: [0,T]\rightrightarrows \Omega$ is monotonically decreasing. 
 \end{compactitem}
 We will investigate the existence of Energetic solutions to the rate-independent systems $(\Spz,\calE_k,\calD)_k$, 
 $k \in \N \cup \{\infty\}$, namely functions $Z: [0,T]\rightrightarrows \Omega$ complying with 
 \begin{compactitem}
 \item[-] the  global stability  condition
\begin{equation}
\label{stab-intro}
\calE_k(t,Z(t))\leq  \calE_k(t,\widetilde Z)+\calD(Z(t),\widetilde Z) \quad \text{for all }
 \widetilde Z \in \Spz, \quad \text{for all } t \in [0,T],
\end{equation}
\item[-] the energy-dissipation balance
\begin{equation}
\label{enbal}
\calE_k(t,Z(t))+\mathrm{Var}_{\calD}(Z;[0,t])=\calE_k(0,Z(0))+\int_0^t\partial_t\calE_k(s,Z(s))\,\mathrm{d}t  \quad \text{for all }  t \in [0,T].
\end{equation}
In \eqref{enbal}, $\mathrm{Var}_{\calD}(Z;[0,t])$  is the total
variation functional induced by the dissipation distance $\calD$, cf.\
\eqref{totvar} below, and the power functional $\partial_t \calE_k$
has to be suitably understood in the brittle case $k=\infty$, in which
the mapping $t\mapsto \calE_\infty(t,Z)$ ceases to be smooth, see
\eqref{power-infty} below. 
\end{compactitem}
As a matter of fact, the  derivative-free character of the Energetic concept makes it suitable to formulate evolutions, like ours,  set up in spaces lacking a linear or even a metric structure, cf.\ e.g.\ the aforementioned pioneering work \cite{DM-Toa2002} on brittle fractures, 
 \cite{MaiMie05EREM} for rate-independent processes in general topological spaces, and \cite{BucButLux} for the rate-independent evolution of debonding membranes. 
\par 
With our first main result, \underline{\bf Theorem \ref{thm:4.2}}, we establish the existence of Energetic solutions for the rate-independent system $(\Spz,\calE_k,\calD)_k$ in the \emph{adhesive case} $k<\infty$. Its proof will be carried out by proving that (the piecewise constant interpolants $(\pwc Z\tau)_\tau$  of) the  discrete solutions $(Z_\tau^i)$ of the associated time-incremental minimization scheme, namely
\begin{equation}
\label{tim-intro}
Z_\tau^i \in \mathop{\mathrm{Argmin}_{Z\in\Spz}}\big\{\calE_k(\sft^i_\tau,Z)+\calD(Z^{i-1}_{\tau},Z)\big\}
\end{equation}
(where $\Pi_\tau:=\{\sft^0_\tau=0<\sft^1_\tau<\ldots,<\sft^N_\tau=T\}$  is  a partition of the time interval $(0,T)$ with step size $\tau$), do converge  to an Energetic solution of the adhesive system as $\tau\down 0$. For this, we will pass to the limit in the discrete versions of the stability condition, and of the upper energy-dissipation estimate $\leq$, to obtain their analogues on the time-continuous level.
\par
The situation for the brittle system is more involved, essentially due  to the intrinsically \emph{nonsmooth} character of the brittle constraint \eqref{brittle-constr-intro}.
This brings about  a nonsmooth time dependence $t\mapsto \calE_\infty(t,Z)$ that has to be handled 
with suitable arguments, since most of the techniques for treating the
Energetic formulation  of rate-independent systems rely on
the condition that the driving energy functional is at least absolutely continuous w.r.t.\ time, cf.\ also Remark \ref{rmk:smooth} ahead.
In our specific case,
because of \eqref{brittle-constr-intro} we are no longer in a position to show that the discrete solutions to \eqref{tim-intro} satisfy a discrete upper energy estimate. Therefore, it remains an \emph{open problem} to prove the existence of Energetic solutions  by passing to the time-continuous limit in scheme \eqref{tim-intro}.  
\par Nonetheless,
time incremental minimization yields the existence of  limiting curves satisfying the stability condition \eqref{stab-intro}. With a terminology borrowed from the theory of gradient flows \cite{AGS08}, we have chosen to qualify such curves  
 as \emph{Stable Minimizing Movements}, cf.\ Definition \ref{def:SMM} ahead. Then, \underline{\bf Theorem \ref{thm:4.1}} asserts the existence of Stable Minimizing Movements both in the adhesive case $k\in \N$ and in the brittle case $k=\infty$.  
\par
The concept of 
Stable Minimizing Movement, though definitely weaker than the Energetic solution notion, seems to be relevant as well. On the one hand, it is 
tightly related to time discretization scheme \eqref{tim-intro}, and it is
on the level of  \eqref{tim-intro} that we can compare our rate-independent evolution with the mean curvature flow, see Remark \ref{ss:2.3}. 
On the other hand, Stable Minimizing Movements enjoy the very same
fine properties as those proved in    \cite{RosTho12ABDM}  for the 
 (semi-)stable delamination variables for the
adhesive contact and brittle delamination systems. Namely, the sets
$Z$ fulfill a \emph{lower density estimate}, which prevents outward
cusps, cf.\ Proposition \ref{prop:LDE}.  Further  geometric properties
of Stable Minimizing Movements are discussed in Section \ref{ss:ulisse}.  
\par
As previously mentioned,  Thm.\  \ref{thm:4.2} will  show that,  
for $k\in \N$ Stable Minimizing Movements
 enhance to Energetic solutions. 
Finally,  \underline{\bf Theorem \ref{thm:4.3}} will provide an existence result for the brittle system. Namely,
we will identify a special setting 
in which the Energetic solutions $(Z_k)_{k}$ of the adhesive systems $(\Spz,\calE_k,\calD)$, 
$k\in\N,$ approximate an Energetic solution of the brittle one $(\Spz,\calE_\infty,\calD)$ as $k\to\infty$.  In this case, we will find that the force term 
$\int_0^T\partial_t\calE_k(t,Z_k(t))\,\mathrm{d}t$ tends to $0$ as $k\to\infty,$ and thus vanishes for the limit system.  
\medskip

\noindent
\paragraph{Plan of the paper.}
In Section \ref{s:2} we fix the setup for the adhesive  and brittle evolution of sets,  and introduce Stable Minimizing Movements and Energetic solutions for both processes. In Section \ref{s:3} we gain further insight into the time incremental minimization scheme \eqref{tim-intro}; for its solutions we  derive a~priori estimates, and discrete versions of the stability condition and of the energy-dissipation inequality.
Section \ref{s:4} contains the statements of all of our existence  results, as well as a detailed discussion 
on the properties of Stable Minimizing Movements,  including some
illustrative numerical experiments in two space
dimensions in the brittle case. Theorems \ref{thm:4.1}, \ref{thm:4.2}, and
\ref{thm:4.3} are then 
proved  in Section \ref{s:6}. 

\medskip

\paragraph{Acknowledgements. }
{\small 
We are very grateful for all the discussions we have shared with Alexander Mielke about rate-independent systems, general mathematical questions, and far beyond. We wish him all the best for the next 60 years with a lot of creative thoughts, and the time to enjoy and elaborate on them. We are already looking forward to be part of this process! 
\par 
Last but not least, the authors also want to acknowledge  financial support: M.T.\ acknowledges the support through the DFG within 
the project ``\emph{Reliability of efficient approximation schemes for material discontinuities described by functions of bounded variation}''
in the priority programme SPP 1748 "\emph{Reliable Simulation
  Techniques in Solid Mechanics. Development of Nonstandard
  Discretisation Methods, Mechanical and Mathematical Analysis}".  R.R.\ has  been
partially supported by  a GNAMPA (INDAM) project. U.S.\  acknowledges
support by the  Vienna Science and Technology Fund (WWTF)
through Project MA14-009 and  by the Austrian Science Fund (FWF)
projects F\,65, P\,27052, and I\,2375.
} 

\section{Setup and notions of solutions}
\label{s:2}
First of all, in the upcoming Section \ref{ss:2.1}  we fix the setting in which the notion of \emph{Energetic solution}  for the adhesive and brittle evolution of sets  can be given. As mentioned in the Introduction, along with Energetic solutions we will also address the much weaker concept of \emph{Stable Minimizing Movement}, which originates from the time-incremental minimization schemes for the adhesive and brittle systems. We shall set up these schemes,  precisely introduce our solution notions, and compare our notion of evolution to the mean curvature flow, in Section \ref{ss:2.2} ahead.  
\subsection{Setup}
\label{ss:2.1}
Throughout this work, $\Omega\subset\R^d$ is a bounded (Lipschitz) domain.
The  evolution of sets  studied in this paper  is driven by a \underline{time-dependent function}
with values in the subsets of $\Omega$, $F: [0,T]\rightrightarrows \Omega$. We will denote by 
$F^c(t)=\Omega\backslash F(t)$ the complement of the set $F(t)$ and impose on the mapping $F$ 
 the following conditions, that also involve a function $f: (0,T)\times\Omega \to \R$, suitably related to $F$:
\begin{subequations}
\label{hypF}
\begin{align}
& 
\label{hypF1}
F(t) \text{ is open for every } t \in [0,T];
\\
& 
\label{hypF2}
F(s) \subset F(t) \quad \text{for every } 0\leq s \leq t\leq T
\\
&
\label{assFsmooth} 
\begin{aligned}
&
\exists\, f \in W^{1,1}(0,T;L^1(\Omega))  \quad \text{such that  }  \foraa\, t\in (0,T)
\\
&
f(t,\cdot) \geq 0 \text{ a.e.\ in }\Omega
\,\text{ and }\,
F(t) = \{ x\in \Omega\, : \ f(t,x)>0 \}\,.
\end{aligned}
\end{align}
\end{subequations}
\begin{remark}
\upshape
\label{rmk:smooth}
Condition \eqref{assFsmooth}  will ensure, for the energy functional
$\calE_k$ for the \emph{adhesive}  system (cf.\ \eqref{en-adh} ahead), that
$\calE_k(\cdot, Z) \in W^{1,1}(0,T)$. This property can be considered
as `standard' within the analysis of rate-independent systems, cf.\
\cite{MieRou-book}. We will rely on  it, for instance,
in the proof of the existence of Energetic solutions in the adhesive case, in order   to readily conclude the energy-dissipation balance once the energy-dissipation upper estimate and the stability conditions have been verified, cf.\ the proof of Thm.\ \ref{thm:4.2} ahead.
\par
 However,  even under \eqref{assFsmooth}, for the brittle system
 the energy functional $\calE_\infty$
 \eqref{en-bri}
  fails to be differentiable w.r.t.\ time, due to the intrinsically nonsmooth character of the brittle constraint \eqref{brittle-constr-intro}.
 One way
to handle a nonsmooth  time-dependence of the driving energy would be to resort to techniques based on the Kurzweil integral, cf.\ \cite{Krejci-Liero} and the references therein.
 In the case of our brittle system, though, we will avoid using
 these  sophisticated tools.  Essentially,  the (somehow) simplified structure of the energy-dissipation balance will allow us to develop ad-hoc arguments in the proof of the existence of Energetic solutions to the brittle system. 
\end{remark}
\par
We now  introduce the  notation for the \underline{state spaces} 
that will enter in our solution concepts for the rate-independent evolution of sets: 
\begin{subequations}
\label{defspaces}
\begin{eqnarray}
\Spz&:=& \{Z\in \mathscr{L}(\Omega)\, : \,\calL^d(Z)<\infty\} =   \mathscr{L}(\Omega)  \,,\\
\Spx&:=&\{ Z\in \mathscr{L}(\Omega)  \, : \,P(Z,\Omega)<\infty\}\,,
\end{eqnarray}
\end{subequations}
where  $\mathscr{L}(\Omega)$ denotes the $\sigma$-algebra of Lebesgue-measurable sets, while   $\calL^d(Z)$ is  the Lebesgue measure of the set $Z$  in $\R^d$, and 
$P(Z,\Omega)$ its perimeter in $\Omega$. 
\par
Since the volume measure $\calL^d$ and the perimeter $P$ are insensitive to null sets, all of our statements will be intended up to null sets. Moreover, 
we will identify  sets  $Z$ from $\Spz$  with the functions $z:\Omega
\to \{0,1\}$ such that 
$\chi_Z = z$, and indeed 
 we will  often use both notations, even within the same line. 
 We recall that 
\begin{equation}
\label{perZ}
\begin{aligned}
P(Z,\Omega)= \sup\left\{ \int_\Omega z \, \mathrm{div}(\varphi)\, : \ \varphi \in \mathrm{C}_{\mathrm{c}}^1(\Omega;\R^d), \ \|\varphi\|_{L^\infty}\leq 1 \right\}  = |\mathrm{D}z|(\Omega),
\end{aligned}
\end{equation}
($\mathrm{C}_{\mathrm{c}}^1(\Omega;\R^d)$ denoting the space of compactly supported $\mathrm{C}^1$-functions on $\Omega$). Hence  
$P(Z,\Omega)<\infty$ if and only if $ z = \chi_Z$ has bounded variation on $\Omega$, i.e.\ $z\in \mathrm{SBV}(\Omega;\{0,1\})$,  since its distributional derivative $\mathrm{D}z$ has no Cantor part. The state spaces for functions  corresponding  to $\Spz$ and $\Spx$ are thus
\begin{eqnarray}
\label{defChispaces}
 L^1(\Omega)\quad\text{ and }\quad
\mathrm{SBV}(\Omega;\{0,1\}):=\{z:\Omega\to\{0,1\} \text{ is the characteristic function of }Z\in\Spx\}\,. 
\end{eqnarray}
Throughout this paper, we will employ  the following notion of convergence of sets: 
we will say that a sequence of sets $(Z_{k})_{k}\subset\Spx$ converges weakly$^{*}$ in $\Spx$ to a limit set 
$Z,$  if for their respective characteristic  functions  $z_k: = \chi_{Z_k}$ and $z : =\chi_Z$, there 
holds weak$^{*}$-convergence in $\SBV(\Omega;\{0,1\}),$ i.e.,
\begin{equation}
\label{sense} 
Z_{k}\conv Z\text{ in }\Spx\quad\Leftrightarrow\quad 
z_k \weaksto z \quad \text{ in } \SBV(\Omega;\{0,1\})\,.
\end{equation}
\par 
Both for the adhesive and the brittle systems we will consider the \underline{dissipation distance}
\begin{align}
\label{dissip}
\calD:\Spz\times\Spz\to[0,\infty],\quad 
\calD(Z_1,Z_2):=\left\{
\begin{array}{ll}
\dissparam \calL^d(Z_1{\setminus} Z_2) &\text{if }Z_2\subset Z_1,\\
\infty&\text{otherwise}
\end{array}
\right.
\end{align}
for a constant $\dissparam>0$. 
For later use, we will also consider the corresponding dissipation distance (and denote it in the same way) on the space of characteristic functions, namely 
\begin{align}
\label{dissip-z}
\calD: 
 L^{1}(\Omega)\times L^{1}(\Omega) \to[0,\infty],\quad 
\calD(z_1,z_2):=\left\{
\begin{array}{ll}
\dissparam(z_1-z_2)&\text{if }z_2\leq z_1 \quad \text{a.e.\ in } \Omega,\\
\infty&\text{otherwise.}
\end{array}
\right.
\end{align}
Clearly, \eqref{dissip}
makes the evolution of a  solution $Z$  to the adhesive/brittle system \emph{unidirectional}, namely $Z: [0,T]\to \Spz$  is nonincreasing:
\begin{equation}
\label{monotonicity}
Z(t)\subset Z(s) 
\quad \text{if }  0\leq s \leq t \leq T.
\end{equation}
In turn,  by this monotonicity property  we have that 
\begin{equation}
\label{totvar}
\begin{aligned}
\mathrm{Var}_{\calD}(Z;[0,t])  & := \sup\left\{ \sum_{j=1}^{M} \calD(Z(\sigma_{j-1}), Z(\sigma_j))\, : \ 0=\sigma_0< \sigma_{1}<\ldots<\sigma_{M-1}<\sigma_M = t\right\} \\\ & = \calD(Z(0),Z(t))\,.
\end{aligned}
\end{equation}
\par
The evolution of the \emph{adhesive} system will be driven by the \underline{energy functional}
\begin{subequations}
\label{ENADH}
\begin{equation}
\label{en-adh}
\calE_k:[0,T]\times\Spx\to[0,\infty)\,,\quad  
\calE_k(t,Z):=P(Z,\Omega)+\calJ_k(t,Z)\quad\text{for }k\in\N\,,
\end{equation}
where 
the functional $\calJ_k:[0,T]\times\Spz \to [0,\infty]$ penalizes the ``brittle constraint'' \eqref{brittle-constr-intro}, namely
\begin{equation}
\label{Jk}
\calJ_k(t,Z):=\int_\Omega k f(t )z\,\mathrm{d}x
\quad \text{for all } (t,Z)\in[0,T]\times \Spz, 
\qquad\text{for }k\in\N\,,
\end{equation}
where $z\in L^1(\Omega)$ is associated with $Z$ via $z=\chi_Z$. 
\end{subequations}
It is immediate to check that $\calE_k$ is differentiable w.r.t.\ time at every $(t,Z) \in [0,T]\times \Spz$, with 
\begin{equation}
\label{power-k}
\partial_t \calE_k(t,Z) = \int_\Omega k \partial_t f(t,x )z(x)\,\mathrm{d}x
\quad \text{for all } (t,Z)\in[0,T]\times \Spz.
\end{equation}
In fact, thanks to \eqref{assFsmooth}  we have that  $\calE_k(\cdot, Z)\in W^{1,1}(0,T)$ for all $Z\in \Spz$. 
With slight abuse of notation, we will sometimes write $\calE_k(t,z)$
(with $P(Z,\Omega)$ rewritten in terms of \eqref{perZ}),
 in place of $\calE_k(t,Z)$.
\par
The energy functional for the \emph{brittle system} is 
\begin{subequations}
\label{ENBRI}
\begin{equation}
\label{en-bri}
\calE_\infty:[0,T]\times\Spx\to[0,\infty]\,,\quad  
\calE_\infty(t,Z):=P(Z,\Omega)+\calJ_\infty(t,Z)
\end{equation}
with  $\calJ_\infty:[0,T]\times\Spz \to [0,\infty]$ the indicator functional associated with  the constraint
\eqref{brittle-constr-intro},  i.e.
\begin{equation}
\label{Jinfty}
\calJ_\infty(t,Z):=
\begin{cases}
0 &\text{ if } Z \cap F(t) =\emptyset,
\\
\infty & \text{ otherwise}.
\end{cases}
\end{equation}
\end{subequations}
Since $F^c(t) = \{ x \in \Omega\, : \, f(t,x) =0 \} $, 
we have that $Z \subset F^c(t)$ if and only if $z = \chi_{Z}$ fulfills
$z(x) f(t,x) =0$ (a.e.\ in $\Omega$). All in all, we have that 
\begin{equation}
\label{J-infty-rephrased}
\calJ_\infty(t,Z) = 
\begin{cases}
0 & \text{ if } z(x) f(t,x)  =0 \quad \foraa\, x \in \Omega,
\\
\infty & \text{ otherwise.}  
\end{cases}
\end{equation}
In fact, also  for the brittle system we will sometimes  write 
$\calJ_\infty$ and 
$\calE_\infty$  as functions of $(t,z)$ through the representation $z=\chi_Z$.   
In place of the usual power functional $\partial_t 
\calE_\infty$,   in this case only the left partial time derivative $\partial_t^- \calE_\infty$ is well defined and fulfills 
\begin{equation}
\label{power-infty} 
\partial_t^-\calE_\infty(t,Z)  :=\lim_{h\uparrow 0}
\frac{\calE_\infty(t+h,Z)-\calE_\infty(t,Z)}{h}  = \lim_{h\uparrow 0}
\frac{\calJ_\infty(t+h,Z)-\calJ_\infty(t,Z)}{h} =0 
\quad \forall\, ( t,Z)\in \mathrm{dom}(\calE_\infty)\,.
\end{equation}
Indeed, $\calE_\infty (t,Z)<\infty$ if and only if $\calJ_\infty(t,Z)=0$, i.e.\ $Z\cap F(t) =\emptyset$. Since $F: [0,T]\rightrightarrows \Omega$ is increasing with respect to time, we then have that $Z\cap F(t+h) =\emptyset$, i.e.\ 
  $\calJ_\infty(t+h,Z)=0$,
  for all $h \in [-t,0)$, which gives \eqref{power-infty}.  Let us stress that the monotonicity property
  \eqref{hypF2} plays a crucial role in ensuring that $\partial_t^-\calE_\infty$ is well defined. 
%
\subsection{Stable Minimizing Movements and Energetic solutions for the adhesive and brittle systems}
\label{ss:2.2}
%
\paragraph{\bf  Time-incremental minimization.} We consider a partition
$\Pi_\tau:=\{\sft^0_\tau=0<\sft^1_\tau<\ldots,<\sft^N_\tau=T\},$ of the time interval $(0,T)$ with step size
$\tau = \max_{i=1,\ldots, N_\tau}(\sft^{i}_\tau {-} \sft_\tau^{i-1})$. 
Discrete solutions arise  from solving  the time-incremental minimization problem:
\emph{starting from $Z_\tau^0: = Z_0 $ for a given $Z_0\in \Spx$, for every $i=1,\ldots, N_\tau$ find}
\begin{equation}
\label{TIM}
Z_\tau^i \in \mathop{\mathrm{Argmin}_{Z\in\Spx}}\big\{\calE_k(\sft^i_\tau,Z)+\calD(Z^{i-1}_{\tau},Z)\big\}\,.
\end{equation}
Here, $k\in \N \cup\{\infty\}$ is fixed and, for simplicity, we choose to omit the dependence of the discrete solutions
$(Z_\tau^i)_{i=1}^{N_\tau}$ 
 on the parameter $k$.
 Time-incremental problem \eqref{TIM} admits a solution 
$Z_\tau^i$ for every $i=1,\ldots, N_\tau$
by the \emph{Direct Method}. Indeed, 
we inductively suppose that the set of minimizers is nonempty at the previous step $i-1$ and consider
  an infimizing sequence $(Z_m)_m$ for the minimum problem at step $i$, with associated functions $(z_m)_m \subset \mathrm{SBV}(\Omega;\{0,1\})$.
Choosing as a competitor in the minimum problem \eqref{TIM}
$Z=\emptyset$ we find that, for every $m\in \N$,
\begin{equation}
\label{est-emptyset}
\begin{aligned}
 \calE_k(\sft^i_\tau,Z_m)+\calD(Z_\tau^{i-1},Z_m) & \leq
\inf_{Z\in\Spx}\big\{\calE_k(\sft^i_\tau,Z)+\calD(Z^{i-1}_{\tau},Z)\big\} + \eps_m\\ 
 & \leq \calE_k(\sft^i_\tau,\emptyset)+\calD(Z^{i-1}_{\tau},\emptyset) + \eps_m \\ & = \dissparam\calL^d(Z^{i-1}_{\tau})  + \eps_m \stackrel{(1)}{\leq} \dissparam\calL^d(Z_0) + \eps_m 
 \end{aligned}
\end{equation}
with $\eps_m\down 0$ as $m\to\infty$.  Note that, for {\footnotesize (1)} we have used that any minimizer  at the  step $i-1$ fulfills  $Z^{i-1}_{\tau}\subset Z_0$  as imposed by $\calD$.
 Since $\calJ_k\geq 0$ for every $k\in \N\cup\{\infty\}$, 
 it is immediate to deduce from estimate \eqref{est-emptyset} that the sequence $(z_m)_m $ is bounded in $ \mathrm{SBV}(\Omega;\{0,1\})$ and, hence, has a weakly-star limit point $\overline{z}$ in $ \mathrm{SBV}(\Omega;\{0,1\})$.
Since $z_m\weaksto \overline{z}$ in $ \mathrm{SBV}(\Omega;\{0,1\})$ implies $z_m\to  \overline{z}$ in $L^q(\Omega)$ for every $1\leq q <\infty$, we have 
\[
\liminf_{m\to\infty} \left( \calE_k(\sft^i_\tau,Z_m)+\calD(Z_\tau^{i-1},Z_m) \right)
\geq  \calE_k(\sft^i_\tau, \overline{Z}) +\calD(z_\tau^{i-1}, \overline{Z})  \qquad \text{for } k \in \N\cup\{\infty\},
\]
so that 
$ Z^{i}_{\tau}: =  \overline{Z}  \in \mathop{\mathrm{Argmin}_{Z\in\Spx}}\big\{\calE_k(\sft^i_\tau,Z)+\calD(Z^{i-1}_{\tau},Z) \big\}$. 
\par
We denote by $\pwc Z\tau: [0,T]\to \Spx$ and   $\upwc Z\tau: [0,T]\to \Spx$ 
the left-continuous and right-continuous piecewise constant interpolants of the elements  
$(Z_\tau^i)_{i=1}^{N_\tau}$, i.e.\
\[
\pwc Z\tau(t): = Z_\tau^i \quad \text{for } t \in (\sft_\tau^{i-1}, \sft_\tau^i], 
\qquad 
\upwc Z\tau(t): = Z_\tau^{i-1} \quad \text{for } t \in [\sft_\tau^{i-1}, \sft_\tau^i) \quad \text{for } i =1,\ldots, N_\tau\
\]
with $\pwc Z\tau(0): = Z_\tau^0$ and $\upwc Z\tau(T) = Z_\tau^{N_\tau}$, while  $\pwl Z\tau$ is  the piecewise linear interpolant
\[
\pwl Z\tau: [0,T]\to \Spx, \qquad \pwl Z\tau(t): = \frac{t-\sft_\tau^{i-1}}{\tau} Z_\tau^i + 
 \frac{\sft_\tau^{i}-t}{\tau} Z_\tau^{i-1}\quad \text{for } t \in [\sft_\tau^{i-1}, \sft_\tau^i] \quad \text{for } i =1,\ldots, N_\tau. 
\]
We will also work with the (left- and right-continuous) piecewise constant interpolants $\pwc \sft\tau:[0,T]\to [0,T]$ and $\underline{\sft}_\tau :[0,T]\to [0,T]$ associated with the partition $\Pi_\tau$. 
\par 
Prior to introducing our solution concepts, we qualify the curves arising as limit points  (in the sense of \eqref{sense}) 
 of the 
interpolants $(\pwc Z{\tau_k})_k$ by resorting to a standard terminology for gradient flows, cf.\ \cite{Ambrosio95, AGS08}. 
\begin{definition}[Generalized Minimizing Movement] 
Let $k\in \N \cup\{\infty\}$.
We call a curve $Z : [0,\widehat T]\to \Spx$, with $0<\widehat T\leq
T$,  a \emph{Generalized Minimizing Movement} starting from $Z_0\in \Spx$ for the rate-independent system $(\Spz,\calE_k,\calD)$ on the interval $[0,\widehat T]$, and write $Z \in \GMMt 0{\widehat T}{k}$,  if 
$Z(0)=Z_0$ and 
there exists a  sequence $\tau_j \down 0$ as $j\to\infty$ such that
\begin{equation}
\label{charact}
\begin{cases}
\exists\, C>0 \ \forall\, j \in \N \ \forall\, t \in [0,\widehat{T}]\,: \quad  \calE(t,\pwc Z{\tau_j}(t)) \leq C, 
\\
\pwc z{\tau_j}(t) \weaksto z(t) \quad \text{in } \mathrm{SBV}(\Omega;\{0,1\}) \ \text{as } j\to\infty \   \text{  for all } t \in [0,\widehat T],
\end{cases}
\end{equation}
with $z(t)=\chi_{Z(t)}$  and $\pwc z{\tau_j}(t) = \chi_{\pwc Z{\tau_j}(t)}$
for all $t\in [0,\widehat T]$.
\par
If $\widehat T=T$, we will simply write $\GMM k$  in place of $\GMMt 0{\widehat T}k$. 
\end{definition}
Observe that every $Z\in  \GMMt 0{\widehat{T}}{k}$ is nonincreasing, i.e.\ \eqref{monotonicity} holds.
Indeed, it is sufficient to observe 
that, for  $j\in \N$
 fixed, $
\pwc Z{\tau_j}(t) \subset \pwc Z{\tau_j}(s) 
$,
since $\pwc Z{\tau_j}(s)  = Z_{\tau_k}^i$ and $\pwc Z{\tau_j}(t) =
Z_{\tau_k}^\ell$ for some $i\leq \ell \in \{0,\ldots,N_j\}$,    and
thus $ Z_{\tau_k}^\ell\subset  Z_{\tau_k}^i$, as the dissipation
distance fulfills $\calD(Z_{\tau_k}^i, Z_{\tau_k}^\ell)<\infty$.  
\par
We are now in a position to introduce 
 the concept of Stable Minimizing Movement. 
\begin{definition}[Stable Minimizing Movement]
\label{def:SMM}
Let $k\in \N \cup \{\infty\}$. We say that a  curve $Z: [0,\widehat{T}]\to\Spx$, with $0<\widehat T\leq T$, is a \emph{Stable Minimizing Movement} starting from  $Z_0 \in \Spx$ for the rate-independent system $(\Spz,\calE_k,\calD)$  on the interval $[0,\widehat T]$, and write $Z \in \SMMt 0{\widehat T}k$,  
 if 
 \begin{enumerate}
 \item
 $Z \in \GMMt 0{\widehat T}k$;
 \item
$Z$ fulfills the stability condition for all $t\in [0,\widehat{T}]$:
\begin{equation}
\label{stab}
\calE_k(t,Z(t))\leq  \calE_k(t,\widetilde Z)+\calD(Z(t),\widetilde Z) \quad \text{for all } \widetilde Z \in \Spx. 
\end{equation} 
\end{enumerate} 
We will  simply  write $\SMM k$ in place of $\SMMt 0{T}k$.
\end{definition}
Requiring the stability condition at all $t\in [0,T]$ clearly implies that the initial datum $Z_0$ will have to be stable at $t=0$, cf.\ \eqref{stab0} ahead. 
\par
In the brittle case $k=\infty$, let us straightforwardly derive  from the stability condition for $k=\infty$ a result on the life-time of Stable Minimizing Movements.
\begin{lemma}
\label{l:elementary}
Suppose that 
$F:[0,T]\rightrightarrows \Omega $ is constant on 
some interval $ [t_1,t_2]\subset [0,T]$.
Let   $Z\in \SMMt 0{t_1}\infty$ fulfill  $\calL^d(Z(t_1))>0$. 
Then, 
the curve $\widehat{Z}: [0,t_2]\rightrightarrows \Omega$ defined by
\[
\widehat{Z}(t) : = \left\{
\begin{array}{ll}
Z(t) &\text{for } t \in [0,t_1],
\\
Z(t_1) &\text{for } t \in (t_1,t_2]
\end{array}
\right.
\]
is  in $\SMMt 0{t_2}\infty$, with 
 $\calL^d(\widehat{Z}(t))>0$ for all $[0,t_2]$. 
\end{lemma}
\par
We now provide a unified definition of Energetic solutions for the adhesive and brittle systems,
where the energy-dissipation balance \eqref{enbal-infty} features the left derivative $\partial_t^- \calE_k$ both for $k \in \N$ and $k =\infty$.  Indeed, in the latter case only $\partial_t^- \calE_\infty$ exists. In the former case, it 
is sufficient to observe that, since for every $k\in \N$  and $Z\in \Spz$ we have that   $\calE_k(\cdot, Z) \in W^{1,1}(0,T)$, there holds
\begin{equation}
\label{used-later}
\partial_t^- \calE_k(t,Z) = \lim_{h\uparrow 0}
\frac{\calE_k(t+h,Z)-\calE_k(t,Z)}{h}   = \partial_t \calE_k(t,Z) \qquad \foraa\, t \in (0,T).
\end{equation}
\begin{definition}[Energetic solution]
Let $k\in \N \cup \{\infty\}$. 
We say that a  curve $Z: [0,T]\to\Spx$ is an \emph{Energetic solution} for the rate-independent system $(\Spz,\calE_k,\calD)$ if it satisfies
\begin{enumerate}
\item the monotonicity property \eqref{monotonicity};
\item   the stability condition \eqref{stab} for all $t\in [0,T]$;
\item the \NEW following energy-dissipation balance 
for all $t\in [0,T]$
\begin{equation}
\label{enbal-infty}
\calE_k(t,Z(t))+\calD(Z(0),Z(t))=\calE_k(0,Z(0))+\int_0^t\partial^-_t\calE_k(s,Z(s))\,\mathrm{d}s \,.
\end{equation} 
\end{enumerate}
\end{definition}
\begin{remark}[Comparison with mean curvature flows]
\label{ss:2.3}
\upshape
Here we  point out some differences of our evolution of sets,   
in  the \emph{adhesive case}, 
to the  (classical) evolution of sets by a mean-curvature flow. 
We perform this comparison in terms of their respective time-discrete schemes. Recall that  scheme \eqref{TIM}, from which the discrete solutions to our adhesive system originate, 
 takes the form
\begin{subequations}
\begin{eqnarray}
\label{ourscheme}
&&Z_{N}^{i}\in\mathrm{argmin}_{\widetilde Z\in\Spz}\left\{P(\widetilde Z,\Omega)
+\int_{\widetilde Z}kf\,\mathrm{d}x+\calD(Z_N^{i-1},\widetilde Z)\right\}
\end{eqnarray}
with $f= f(\mathsf{t}_\tau^i)$. 
We highlight that, here, due to  
  \eqref{dissip} the dissipation potential 
  accounts for unidirectionality of the
evolution
by enforcing that $Z_{N}^{i} \subset Z_{N}^{i-1}$. 
 This is a first difference to classical mean curvature flows, which do not take into account a unidirectional evolution. 
 \par
 Thus, for further comparison let us for the moment 
disregard  unidirectionality and confine the discussion to a \emph{symmetric} dissipation distance, i.e.\ 
\begin{equation}
\label{D-symmetric}
\calD(\widetilde Z,Z)=\calD(Z,\widetilde Z)=\int_{\Omega} \dissparam |\tilde z-z|\,\mathrm{d}x\,.
\end{equation} 
Since $\tilde z,z$ are the characteristic functions of the finite-perimeter sets $\widetilde Z,Z,$ and thus only take the values $0$ or $1,$ we can equivalently
 rewrite our dissipation distance as a \emph{squared} distance, i.e.\ 
 $\calD(\widetilde Z,Z)=\calD(Z,\widetilde
 Z)=\int_{\Omega}\dissparam|\tilde z-z|^{2}\,\mathrm{d}x$. This 
 makes our scheme \eqref{ourscheme} closer to that for the mean
 curvature flow, which is usually related to a quadratic, symmetric  distance.  More precisely,  
 with \eqref{D-symmetric}  our time-discrete scheme rephrases as 
\begin{eqnarray}
\label{general}
Z_{N}^{i}\in\mathrm{argmin}_{\widetilde Z\in\Spz}\left\{P(\widetilde Z,\Omega)
+\int_{\widetilde Z}kf\,\mathrm{d}x
+\displaystyle\frac{N}{T}\int_{\Omega}\dissparam|\chi_{Z_N^{i-1}}-\chi_{\widetilde Z}|^{2}\,\mathrm{d}x\right\}
\end{eqnarray}
\par
Depending on the values of $f$, minimization problem \eqref{general} may allow for an infinite number of minimizers. This is e.g.\ the case if $f=\mathrm{const}$ in a large open connected set of positive measure,  because a minimizer that is translated by a sufficiently small distance is still a minimizer of \eqref{general}. 
To make a selection of minimizers that keeps the minimizer pinned, following the classical literature on mean curvature flows, cf.\ e.g.\ \cite{ATW93CDFV,LuckStur95ITDM,Visi97MMCN,Visi98NMCF}, one rather replaces the above quadratic dissipation distance by the following expression $\int_{\Omega}\alpha(-\tfrac{N}{T}\,\mathrm{sdist}(x,\partial Z_N^{i-1}))\tilde Z\,\mathrm{d}x,$ 
where $\mathrm{sdist}(x,\partial E)=\mathrm{ess\,inf}\{|x-y|,\,y\in\Omega\backslash E\}
-\mathrm{ess\,inf}\{|x-y|,\,y\in E\}$ denotes the signed distance. The
case $\alpha:\R\to\R$ nonconstant, bounded, and monotone is discussed in the works \cite{Visi97MMCN,Visi98NMCF}, while 
$\alpha:\R\to\R$ linear is the original and well-established ansatz first proposed in \cite{ATW93CDFV}. 
We now combine this specific discrete dissipation distance, multiplied with a prefactor $\eps>0,$ 
with our choice to obtain the discrete problem 
\begin{eqnarray}
\label{Visintin}
Z_{N}^{i}\in\mathrm{argmin}_{\widetilde Z\in \Spz}&\bigg\{P(\widetilde Z,\Omega)
+\!\! \displaystyle\int_{\widetilde Z}kf\,\mathrm{d}x
+\!\!\displaystyle\int_{\Omega}\!\!\dissparam |\chi_{Z_N^{i-1}}-\chi_{\widetilde Z}|^2 \,\mathrm{d}x
\nonumber\\
&\qquad \quad +\eps \displaystyle\!\!\int_{\Omega}\!\!\alpha\left(-\displaystyle\frac{N}{T}\,\mathrm{sdist}(x,\partial Z_N^{i-1})\right)\tilde Z\,\mathrm{d}x\bigg\}\,.\;
\end{eqnarray}
\end{subequations}
For fixed $\eps>0$, this minimization problem is  a particular case of that addressed  in  \cite{Visi98NMCF}. 
Therefore,  our time-discrete problem with the symmetric dissipation distance  from \eqref{D-symmetric}, i.e.\ \eqref{general},   corresponds to the limit case 
$\eps=0$  (so that the term with the signed distance function disappears), in \eqref{Visintin}.
 Thus, formally, our scheme \eqref{general} can be understood as a singularly perturbed limit of the flow \eqref{Visintin} set forth in  
 \cite{ATW93CDFV, Visi98NMCF}. 
\end{remark}

%
\section{The time-discrete problem}
\label{s:3}
In this section we
show that the time-discrete solutions arising from scheme \eqref{TIM} satisfy approximate versions of the stability condition and of the lower energy-dissipation estimates. Instead, as we will see, due to the
unidirectionality constraint on the evolution  we will be able to
obtain only a \emph{discrete} energy-dissipation estimate 
under the restriction that $k\in\N$, i.e.\ for  adhesive systems. 
%
\par
By testing the minimality of $Z^i_{\tau}$ at time step $\sft^i_\tau$ (cf.\ scheme \eqref{TIM}), with any $\widetilde Z\in\Spx,$ and by exploiting
that the dissipaton distance $\calD$ satisfies the triangle inequality, i.e.\
\begin{align*}
\calE_k(\sft^i_\tau,Z^i_{\tau})+\calD(Z^{i-1}_{\tau},Z^i_{\tau})
\leq\calE_k(\sft^i_\tau,\widetilde Z)+\calD(Z^{i-1}_{\tau},\widetilde Z)
\leq\calE_k(\sft^i_\tau,\widetilde Z)+\calD(Z^i_{\tau},\widetilde Z)+\calD(Z^{i-1}_{\tau},Z^i_{\tau})\,,
\end{align*}
we can show that the time-incremental solutions satisfy the following {\bf stability condition}:
\begin{equation}
\label{discrstab}
\forall\,\widetilde Z\in\Spx:\quad
\calE_k(\sft^i_\tau,Z^i_{\tau})\leq \calE_k(\sft^i_\tau,\widetilde Z)+\calD(Z^i_{\tau},\widetilde Z)\,.
\end{equation}
\par 
From \eqref{discrstab}, choosing $\widetilde Z=\emptyset$ as a competitor and arguing as for \eqref{est-emptyset}, we deduce  
the following uniform a priori bound for the time-incremental
minimizers:
\begin{equation}
\label{unibdstab}
\calE_k(\sft^i_\tau,Z^i_{\tau})\leq\calE_k(\sft^i_\tau,\emptyset)+\calD(Z^i_{\tau},\emptyset)\leq\calD(Z_0,\emptyset)=\dissparam \calL^d(Z_0)\,.
\end{equation}
\par
Moreover, testing stability at time $\sft_{i-1}^\tau$ with $Z_i^{\tau}$ we obtain a {\bf discrete
lower energy estimate}:
\begin{equation}
\label{discr-lee}
\begin{split}
\calE_k(\sft^{i-1}_\tau,Z^{i-1}_{\tau})
&\leq\calE_k(\sft^{i-1}_{\tau},Z^i_{\tau})+\calD(Z^{i-1}_{\tau},Z^i_{\tau})\\
&\stackrel{(2)}{=}\left\{
\begin{array}{ll}
\calE_\infty(\sft^i_\tau,Z^i_{\tau})+\calD(Z^{i-1}_{\tau},Z^i_{\tau})&\text{for }k=\infty\,,\\
\calE_k(\sft^i_{\tau},Z^i_{\tau})+\calD(Z^{i-1}_{\tau},Z^i_{\tau})-\int_{\sft^{i-1}_{\tau}}^{\sft^i_\tau}\partial_t\calE_k(t,Z^i_{\tau})\,\mathrm{d}t
&\text{for }k\in\N\,.          
\end{array}\right.
\end{split}
\end{equation}
In the brittle case  $k=\infty,$
equality {\footnotesize (2)} 
 is due to the fact that $\calJ_\infty(Z^i_\tau,F(\sft^{i-1}_\tau))=\calJ_\infty(Z^i_\tau,F(\sft^i_\tau))=0$ since 
  $F(\sft^{i-1}_\tau)\subset F(\sft^{i}_\tau)$ by monotonicity of $F$. 
 Instead, in the adhesive case  $k\in\N$
  we indeed have 
\begin{equation}
\calE_k(\sft^{i-1}_\tau,Z^i_{\tau})-\calE_k(\sft^{i}_\tau,Z^i_{\tau})
=-\int_{\sft^{i-1}_\tau}^{\sft^i_\tau}\partial_t\calE_k(t,Z^i_{\tau})\,\mathrm{d}t\,.
\end{equation}
All in all, taking into account \eqref{used-later} and the fact 
 that 
$\partial_{t}^{-}\calE_{\infty}(t,Z^{i}_{\tau})=0$ by \eqref{power-infty},
estimate \eqref{discr-lee} can be 
rewritten as
\begin{equation}
\calE_k(\sft^{i-1}_\tau,Z^{i-1}_{\tau})\leq \calE_k(\sft^i_{\tau},Z^i_{\tau})+\calD(Z^{i-1}_{\tau},Z^i_{\tau})-\int_{\sft^{i-1}_{\tau}}^{\sft^i_\tau}\partial_t^{-}\calE_k(t,Z^i_{\tau})\,\mathrm{d}t \qquad \text{for all }  k\in\N\cup\{\infty\}\,.
\end{equation}
\par
Finally, in the adhesive case $k\in\N$
 one can also obtain a {\bf discrete upper energy-dissipation estimate}
by testing the minimality of $Z^i_{\tau}$ at time $\sft^i_\tau$ by $\widetilde Z: = Z^{i-1}_{\tau}$, i.e.\ we get
\begin{equation}
\label{uedek}
\calE_k(\sft^i_{\tau},Z^i_{\tau})+\calD(Z^{i-1}_{\tau},Z^i_{\tau})\leq \calE_k(\sft^i_\tau,Z^{i-1}_{\tau})= \calE_k(\sft^{i-1}_\tau,Z^{i-1}_{\tau})
+\int_{\sft^{i-1}_{\tau}}^{\sft^i_\tau}\partial_t\calE_k(t,Z^{i-1}_{\tau})\,\mathrm{d}t\quad\text{for }k\in\N\,.
\end{equation}
Instead, for $k=\infty$ the constraint 
$F(\sft^i_\tau) \cap \widetilde Z =\emptyset$ imposed by the 
  energy contribution 
$\calJ_\infty(\sft_\tau^i,\cdot)$ 
forbids us to choose $\widetilde Z: = Z^{i-1}_{\tau}$ \NEW as a competitor in the minimization problem \eqref{TIM}.  
In fact,  we have 
$F(\sft^{i-1}_\tau)\subset F(\sft^i_\tau)$ and thus 
$ Z^{i-1}_{\tau} \subset F^c(\sft^{i-1}_\tau)$ need not satisfy the constraint  $ Z^{i-1}_{\tau} \subset F^c(\sft^{i}_\tau)$.   
\par 
We are now in a position to deduce from the   above observations  
(summing up the discrete lower and upper energy inequalities \eqref{discr-lee} and \eqref{uedek} over the index $i$), the following result. 
\begin{proposition} 
\label{DiscrProps}
Consider the rate-independent systems $(\Spz,\calE_k,\calD)$ for $k\in\N\cup\{\infty\}$ as defined 
by \emph{\eqref{hypF}}, \emph{\eqref{dissip}}, and \emph{\eqref{ENADH}} if $k\in \N$, \emph{\eqref{ENBRI}} if $k=\infty$.
Then,
\begin{compactenum}
\item  
For every $k\in\N\cup\{\infty\}$ and every $t\in(0,T]$, 
the interpolants 
$(\pwc Z\tau)_\tau$
satisfy the time-discrete 
stability condition 
\begin{equation}
\label{discr-stab-interp}
\calE_k(\pwc \sft\tau(t),\pwc Z\tau(t)) \leq \calE_k(\pwc \sft\tau(t),\widetilde Z) +\calD(\pwc Z\tau(t),\widetilde Z) \quad \text{for all } \widetilde Z \in \Spz, 
\end{equation}
 as well as  the uniform a priori bound 
 \begin{equation}
 \label{energy-bound-interp}
 \exists\, C>0  \ \forall\, k \in \N\cup\{\infty\} \  \forall\,\tau>0\, : \qquad 
 \sup_{t\in (0,T)}\calE_k(\pwc \sft\tau(t),\pwc Z\tau(t)) \leq C\,.
 \end{equation}
\item 
For every $k\in\N\cup\{\infty\} $  and every $t\in[0,T]$
there holds the  \emph{lower} energy-dissipation estimate
\begin{equation}
\label{lower-discrete}
\begin{split}
\calE_k(0,Z_0)
+\int_{0}^{\pwc \sft\tau(t)}\partial_t^{-}\calE_k(r,\pwc Z{\tau}(r))\,\mathrm{d}r 
\leq \calE_k(\pwc \sft\tau(t), \pwc Z\tau(t))+\calD(Z_0,\pwc Z\tau(t))\,.
\end{split}
\end{equation}
\item For every   $k\in\N$
and every $t\in [0,T]$ there holds the 
\emph{two-sided} energy-dissipation estimate 
\begin{equation}
\label{enbd}
\begin{split}
\calE_k(0,Z_0)
+\int_{0}^{\pwc \sft\tau(t)}\partial_t\calE_k(r,\pwc Z{\tau}(r))\,\mathrm{d}r 
 & \leq
\calE_k(\pwc \sft\tau(t), \pwc Z\tau(t))+\calD(Z_0,\pwc Z\tau(t))
\\
&\leq
\calE_k(0,Z_0)
+\int_{0}^{\pwc \sft\tau(t)}\partial_t\calE_k(r,\upwc Z{\tau}(r))\,\mathrm{d}r .
\end{split}
\end{equation}
\end{compactenum} 
\end{proposition}

\section{Main results}
\label{s:4}
Our first result ensures that, both for the adhesive and the brittle systems,
the set $\GMM k \neq \emptyset$, and that  Generalized Minimizing Movements starting from stable initial data are in fact Stable Minimizing Movements. 
\begin{theorem}
\label{thm:4.1}
Let $k\in\N\cup\{\infty\}$.
 Let the rate-independent systems $(\Spz,\calE_k,\calD)$  fulfill
 \emph{\eqref{hypF}}, \emph{\eqref{dissip}}, and \emph{\eqref{ENADH}}
 if $k\in \N$, \emph{\eqref{ENBRI}} if $k=\infty$.
Then,
\begin{enumerate}
\item $\GMM k \neq \emptyset$ 
 for every $Z_0\in \Spx$, 
and for every $Z \in   \GMM k $ and for any sequence $(\pwc z{\tau_j})_j$ fulfilling 
\eqref{charact}
 there also holds
 \begin{subequations}
 \label{refined-convs}
 \begin{align}
 &
  \label{refined-convs-1}
 \pwc z{\tau_j}(t), \,  \upwc z{\tau_j}(t)  \weaksto z(t)  &&   \text{ in } \SBV(\Omega;\{0,1\}) \text{ for almost all } t \in (0,T);
 \\
 &
  \label{refined-convs-2} 
 \pwc z{\tau_j}, \, \upwc z{\tau_j}  \weaksto z  &&  \text{ in } L^\infty (0,T;\SBV(\Omega;\{0,1\}),  
 \\
 &
  \label{refined-convs-3} 
 \pwc z{\tau_j}, \, \upwc z{\tau_j}  \to z  &&  \text{ in } L^q ((0,T){\times} \Omega) 
\text{ for every } q \in [1,\infty).
\end{align}
 \end{subequations}
\item If  in addition
\begin{equation}
\label{stab0}
\calE_k(0,Z_0) \leq \calE_k(0,\widetilde Z) +\calD(Z_0,\widetilde Z) \quad \text{for all } \widetilde Z \in \Spx,
\end{equation}
then every Generalized Minimizing Movement for the rate-independent system $(\Spz,\calE_k,\calD)$   is also a  Stable Minimizing Movement,  i.e.\ $\GMM k =\SMM k$. 
\end{enumerate}
\end{theorem}
We postpone the \emph{proof} of Theorem \ref{thm:4.1} to Section \ref{s:6}.  Proposition \ref{prop:LDE} below  ensures that, both in the adhesive and in the brittle cases,  any Stable Minimizing Movement $Z$
 enjoys a regularity property,
 which prevents the occurrence of outward cusps,
  introduced by \textsc{S.\ Campanato} as  Property $\frak{a},$ 
cf.\ e.g.\ \cite{Camp63PHAC,Camp64PUFS},
and also  known as  \emph{lower density estimate} in e.g.\ \cite{FonFra95RBVQ,AmFuPa05FBVF}. We recall it in the following definition. 
\begin{definition}[Property $\frak{a}$] 
\label{def-propa}
A set $M\subset\R^d$ has Property $\frak{a}$ if there exists a constant $\frak{a}>0$ such that
\begin{align}
\label{propa}
\forall\,y\in M\;\;\forall\,\rho_\star>0:\quad 
\calL^{d}(M\cap B_{\rho_\star}(y))\ge
 \frak{a} \rho_\star^{d}\,.
\end{align}
\end{definition}
\NEW In \cite{RosTho12ABDM} we were able to prove the validity of Property $\frak{a}$
for Energetic solutions to  adhesive contact  and brittle delamination systems, in which  $z = \chi_Z$ is a phase-field parameter describing the state of the bonds between two bodies.  
This analysis can be extended to the more general context of
our evolution of sets. 
More precisely, 
 we can show that  
any set $Z\in \SMM k$ 
 fulfills \eqref{propa}, with constants \emph{uniform} w.r.t.\ the parameter $k$, cf.\ \eqref{LDE} below, and at all points in  the support of its characteristic function $z=\chi_Z$, 
defined in measure-theoretic way as 
\begin{equation}
\label{defsupp}
   \supp z :=\bigcap\{A \subset 
\R^{d} \, : \,A\text{ closed },\, \calL^{d} (Z{\backslash} A)=0\}.
   \end{equation}
However,  for this we need to additionally impose that $\Omega$ is \emph{convex}: this is essential for the proof of the  uniform relative isoperimetric inequality from \cite[Thm.\ 3.2]{Thom15}, which is in turn a key ingredient for obtaining \eqref{LDE}, cf.\ 
\cite[Sec.\ 6]{RosTho12ABDM}. 
\begin{proposition}
\label{prop:LDE}
Let $k\in\N\cup\{\infty\}$ and  the rate-independent systems $(\Spz,\calE_k,\calD)$ fulfill
 \emph{\eqref{hypF}}, \emph{\eqref{dissip}}, and \emph{\eqref{ENADH}} if $k\in \N$, \emph{\eqref{ENBRI}} if $k=\infty$.
Suppose in addition that 
\begin{equation}
\label{Omega-cvx}
\Omega \text{ is convex.}
\end{equation}
 Then, every 
 $Z\in \SMM k$  satisfies the following  lower density estimate:
there are constants $R$ and $\mathfrak{a} = \frak{a}(\Omega,d, \dissparam)>0$ 
 depending solely on $\Omega\subset\R^{d},$ space dimension $d,$ 
and on the parameter $\dissparam,$  
such that   for every $k\in \N \cup \{\infty\}$ there holds
\begin{align}
\label{LDE}
\forall\,y\in\supp z \quad \forall\,\rho_\star>0:\qquad 
\calL^{d}(Z\cap B_{\rho_\star}(y))\ge
\begin{cases} 
\mathfrak{a} \rho_\star^{d}&\text{if }\rho_\star< R,\\
\mathfrak{a}  R^{d}&\text{if }\rho_\star\geq R\,.
\end{cases}
\end{align}
\end{proposition}
The proof of Proposition \ref{prop:LDE} follows from directly  adapting the argument
developed in  \cite[Sec.\ 6]{RosTho12ABDM}, to which we refer the reader. 
\par
A straightforward consequence of Proposition \ref{prop:LDE} is Corollary \ref{cor:k} below, stating  that, the first time $t_*$ at which the complement set $F^c(t_*)$ violates the  volume constraint of the  lower density estimate \eqref{LDE} is an extinction time for Stable Minimizing Movements in the brittle case.  
\begin{corollary}
\label{cor:k}
Let $k=\infty$.  Let the rate-independent system $(\Spz,\calE_\infty,\calD)$ fulfill
 (\ref{hypF}), (\ref{dissip}), and (\ref{ENBRI}) if $k=\infty$.
Additionally, assume \eqref{Omega-cvx}.
  Suppose that  there is  
 $t_* \in (0,T]$ such that 
\begin{equation}
\label{F-tstar-c}
\calL^d(F^c(t_*))<\frak{a}(\Omega)R^d
\end{equation}
 for $R$ from \eqref{LDE}. 
Then, every $Z\in \SMM \infty$ fulfills   $\calL^d(Z(t_*))=0$, and 
consequently  $\calL^d(Z(t))=0$ for all $t>t_*$. 
\end{corollary}
Indeed, from $Z(t_*) \subset F^c(t_*)$ we gather that $\calL^d(Z(t_*))<\frak{a}(\Omega)R^d$, therefore  $Z(t_*)$  violates \eqref{LDE}. Then, $ \calL^d(Z(t_*))=0$, which implies that  $\calL^d(Z(t))=0$ for all $t>t_*$ due to the monotonicity property $Z(t)\subset Z(t_*)$. Observe that this argument strongly relies on the brittle constraint \eqref{brittle-constr-intro}; in fact, it is not clear how to obtain an analogue of Corollary \ref{cor:k}  for the adhesive system.
We shall  provide a more detailed discussion of further properties of Stable Minimizing Movements in Section \ref{ss:ulisse}. 
\par
In the adhesive case $k\in \N$, Stable Minimizing Movements enhance to Energetic solutions under the very same conditions as for the existence Thm.\ \ref{thm:4.1}. 
\begin{theorem}
\label{thm:4.2}
Let $k\in \N$.   Let the rate-independent systems $(\Spz,\calE_k,\calD)$  fulfill
 (\ref{hypF}), (\ref{dissip}), and (\ref{ENADH}), and let $Z_0\in \Spx$ comply with the stability condition at $t=0$, cf.\ \eqref{stab0}.  
 Then, every $Z\in \SMM k $ is    an Energetic solution to the rate-independent system $(\Spz,\calE_k,\calD)$.  
\end{theorem}
The proof of Thm.\ \ref{thm:4.2}, postponed to Section \ref{s:6}, will be carried out by  passing to the limit in the time-discrete scheme  \eqref{TIM} in the following steps:  it will be shown that any $Z\in \SMM k $ complies with the upper energy-dissipation estimate by passing to the limit in its discrete version \eqref{enbd} via lower semicontinuity arguments; the lower energy estimate will then follow from the stability condition, via a by-now standard technique. 
\par
We are not in a position to prove the analogue of Thm.\ \ref{thm:4.2} in the  brittle case 
$k=\infty$,  since the discrete version of the upper energy-dissipation estimate is not at our disposal, and therefore   the upper 
estimate in the energy-dissipation balance \eqref{enbal-infty} cannot be obtained by passing to the limit in the time-discretization scheme.  The existence of Energetic solutions to the brittle system can be proven, though, by  taking the limit as $k\to\infty $ in the Energetic formulation at the time-continuous level, provided that the `external force'  $F: [0,T]\rightrightarrows \Omega$ additionally fulfills  condition
\eqref{power-control} below, and that the initial datum $Z_0$ complies with a suitable compatibility condition, cf.\ \eqref{compatibility} below. This is stated in Theorem \ref{thm:4.3} ahead where, 
for completeness, we also give a result on the convergence of Stable Minimizing Movements
for the adhesive system to 
Stable Minimizing Movements for the brittle one. The proof of Thm.\ \ref{thm:4.3} shall be also performed in Sec.\ \ref{s:6}.  
\begin{theorem} 
\label{thm:4.3}
Let the initial data of the adhesive problems $Z_0^{k}$ fulfill \eqref{stab0} for each $k\in\N$ and assume that 
\begin{equation}
\label{well-prep1}
Z_{0}^{k}\conv Z_{0}\text{ in }\Spx\quad\text{in the sense of \eqref{sense}}\,.
\end{equation}
Then,
\begin{enumerate}
\item any sequence $(Z_k)_k$ with $Z_k \in \SMMk$ for every $k\in \N$ admits a (not relabeled) subsequence 
such that $Z_{k}(t)\conv Z(t)$ in $\Spx$ in the sense of \eqref{sense}  for all $t\in [0,T]$.   
\item 
Assume that 
\begin{equation}
\label{power-control}
\partial_t f\geq0  \qquad \text{a.e.\ in  } [0,T]\times\Omega,
\end{equation}
and  that the initial data $(Z_{0}^{k})_{k}$ of the adhesive problems are well-preprared, i.e., in addition to \eqref{well-prep1} there holds 
\begin{equation}
\label{well-prep}
\calE_k(0,Z_0^k) \to \calE_\infty (0,Z_0) \qquad \text{as } k \to\infty. 
\end{equation}
Further, assume that the limit initial datum $Z_{0}$ satisfies the compatibility condition 
\begin{equation}
\label{compatibility}
P(Z_{0},\Omega)=\calD(Z_{0},\emptyset) = \dissparam \calL^d(Z_0) \,. 
\end{equation}
Then any sequence $(Z_k)_k$ of Energetic solutions of the adhesive systems $(\Spz,\calE_k,\calD)$ 
with $Z_{k}(0)=Z_{0}$ admits a (not relabeled) subsequence converging  as $k\to\infty$, in the sense of \eqref{sense}, to an Energetic solution of the brittle system $(\Spz,\calE_\infty,\calD),$ i.e.\ energy-dissipation balance \eqref{enbal-infty} holds true.   Moreover, the analogue  of the compatibility condition \eqref{compatibility}  holds true also 
for all $t\in(0,T],$ i.e.,
\begin{equation}
\label{compatibility-t}
P(Z(t),\Omega)=\calD(Z(t),\emptyset)= \dissparam \calL^d(Z(t))\,. 
\end{equation}
\end{enumerate}
\end{theorem}
A few comments on the above statement are in order:
\begin{enumerate}
\item A close perusal of the proof of Lemma \ref{PropGamma} ahead on the $\Gamma$-convergence of the functionals $(\calE_k)_k$ to $\calE_\infty$  reveals that, under \eqref{well-prep1} there holds
$ P(Z_{0},\Omega) \leq \liminf_{k\to\infty} P(Z_{0}^{k},\Omega)$ and $\calJ_\infty (0,Z_0) \leq  \liminf_{k\to\infty} \calJ_k(0,Z_0^k)$. Therefore, \eqref{well-prep} is indeed equivalent to requiring that
\begin{equation}
\label{conseq-well-prep}
P(Z_{0}^{k},\Omega)\to P(Z_{0},\Omega) \quad \text{ and } \quad    \calJ_k(0,Z_0^k) \to \calJ_\infty (0,Z_0) =0 \qquad \text{as } k \to \infty\,.
\end{equation}
\item 
Observe that  the `brittle energy-dissipation balance'  in fact reads
\begin{equation}
\label{proprio-cosi}
P(Z(t),\Omega) + \dissparam \calL^d(Z_0{\setminus}Z(t)) = P(Z_0,\Omega) \qquad \text{for all } t \in (0,T]
\end{equation}
taking into account that $\calJ_\infty(t,Z(t)) = \calJ_\infty(0,Z_0)=0$ and that $\int_0^t \partial_t^-\calE_\infty(s,Z(s)) \dd s =0$. Combining \eqref{proprio-cosi} with \eqref{compatibility}, we immediately conclude \eqref{compatibility-t}. 
The compatibility condition \eqref{compatibility} 
 may lead to earlier extinction as  Example 
 \ref{ex:4.7} below illustrates. Examples of nontrivial 
 evolutions complying with \eqref{compatibility} at all times $t\in [0,T]$ are provided with Examples \ref{ex:4.8} ahead. 
 Finally, we also comment on the outcome of condition
\eqref{compatibility} for the adhesive system later on in  Remark \ref{rmk:compat4adh}.  
\end{enumerate}
\begin{example}[Extinction of sets under compatibility condition \eqref{compatibility}] 
\label{ex:4.7}
\upshape
Let $Z_{0}\subset\Omega$ be  a ball of radius $r$ in $\R^{2}$. Compatibility condition 
\eqref{compatibility} holds true if $\calD(Z_{0},\emptyset)=\dissparam\pi r^{2}=P(Z_{0},\Omega)=2\pi r$. 
This is satisfied for $r=2/\dissparam$ and for $r=0$. In other words, an Energetic solution can only exist for 
the initial ball $Z_{0}$ with $r=2/\dissparam$.  As soon as the ball is forced by the brittle constraint to shrink, it is extinguished. 
 In terms of Stable Minimizing Movements, the ball $Z(t)$ would be extinguished according to \eqref{F-tstar-c} 
only if $\calL^{2}(Z(t))<\frak{a}(\Omega)R^{2},$ which may be a
smaller value than enforced  by the compatibility condition \eqref{compatibility}.  
\par 
Yet, the above considerations can be used to provide an example of an Energetic solution for a rate-independent evolution of sets in $\Omega\subset\R^{2}$: 
Consider $m$ balls $B_{2/\dissparam}(x_{i}),$ $i=1,\ldots,m,$ with centers $x_{i}\in\Omega$ such that the balls are pairwise disjoint and such that 
$\cup_{i=1}^{m}B_{2/\dissparam}(x_{i})\subset\Omega$. Along the time interval $[0,T]$ we choose a partition $t_{0}=0<t_{1}<\ldots<t_{m}=T$ and we define the evolution of the forcing set $F$ through its complement $F^{c}$ by setting $F^{c}(t):=\cup_{i=1}^{m-k}B_{2/\dissparam}(x_{i})$ for all $t\in[t_{k},t_{k+1}),$ for  $k=0,1,\ldots,m-1$. In this case, $F^{c}(t)$ itself is a Stable Minimizing Movement, i.e.\  $Z(t)=F^{c}(t),$ 
which additionally satisfies the compatibility condition \eqref{compatibility} and thus provides an Energetic solution of the system.   
\end{example}

We devote the remainder of this section to a collection of remarks,  illustrating some more 
 properties of Stable Minimizing Movements for the \emph{brittle} system.
%
\subsection{Basic features of Stable Minimizing Movements for the brittle system}
\label{ss:ulisse}
In order to gain further insight into the features of Stable
Minimizing Movements for the brittle system, we shall address a single
step of the Minimizing Movement procedure, starting from some
given initial set $Z^o\subset\Omega$ under the action of the forcing $F$. Along the whole
subsection we assume that $F^c \subset Z^o$. Hence,  the {\it
  single-step} minimization problem 
  \[
  Z\in \text{\rm argmin} \{P(Z, \Omega)+\dissparam\calL^d(Z^o{\setminus}Z) \ : \ 
Z \in \mathbf{X}, \ 
 Z  \subset  F^c \cap Z^o\},
  \]
 can be equivalently reformulated as 
\begin{equation}
Z\in \text{\rm argmin} \{P(Z, \Omega)+\dissparam\calL^d(F^c{\setminus} Z) \ : \ 
Z \in \mathbf{X}, \ 
Z  \subset F^c\},
\label{c1}
\end{equation}
taking into account that $\calL^d(Z^o{\setminus}Z) = \calL^d(Z^o{\setminus}F^c) + \calL^d(F^c{\setminus}Z)$. Observe that 
problem \eqref{c1} has the advantage of being in fact independent of
$Z^o$. 
In the following, we comment on 
  some properties of the minimizers of \eqref{c1}.

\begin{figure}[h]
  \centering
  \pgfdeclareimage[width=55mm]{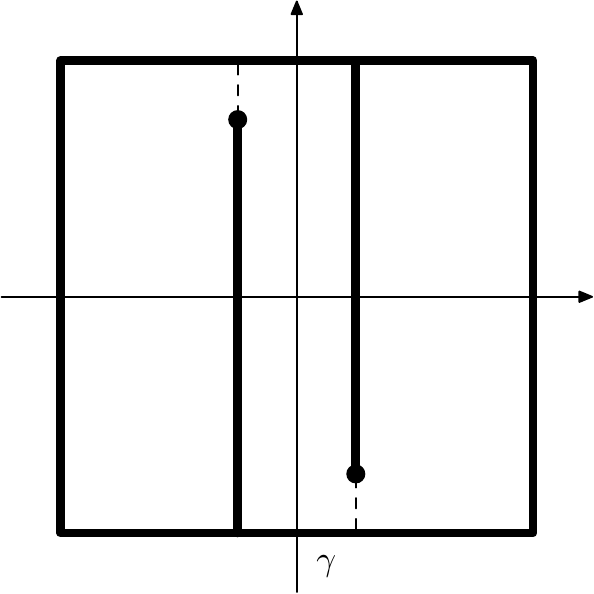}{cut.pdf}
    \pgfuseimage{cut.pdf}
  \caption{ An example for nonconnectedness: A needle-like forcing $F$.}
\label{cut}
\end{figure}

\paragraph{\bf Connectedness.} 
Even starting from a connected initial set
$Z^o$, connectedness is not necessarily preserved 
by Stable Minimizing Movements. Indeed, it is not preserved by solutions  of \eqref{c1}. 
An example in this direction is given by the incremental
minimizer under the effect of the needle-like forcing  $F$ (see Figure \ref{cut}) given by 
\begin{equation} 
F:=\big(\{\gamma\}\times  [-1+\gamma,\infty)\big) \cup \big(\{-\gamma\}\times
 (-\infty,1-\gamma]\big)  \cup ([-1,1]^2)^c \label{fff}
\end{equation}
for some suitably small  $\gamma\in (0,1),$ tuned to $\dissparam$ and
specified later on.
Assume by contradiction that $Z$ is connected
and consider the disconnected competitor $Z^{\mathrm{disc}} :=Z\setminus
([-\gamma,\gamma]\times[-1,1])$. 
We have that 
$$\calL^{2}(F^c {\setminus}Z^{\mathrm{disc}}) = \calL^2( F^c {\setminus}Z)+\calL^2(Z{\setminus}Z^{\mathrm{disc}}) 
\leq \calL^2( F^c {\setminus}Z) + 4\gamma,$$
 where we have estimated $\calL^{2}(Z\backslash
Z^{\mathrm{disc}}) \leq \calL^{2}([-\gamma,\gamma]\times[-1,1])=4\gamma$. 
On the one hand, by passing from $Z$ to $Z^{\mathrm{disc}},$ the perimeter drops at
least  by twice the distance between the points
$(-\gamma,1-\gamma)$ and $(\gamma,-1+\gamma),$ which is
$2\sqrt{(2\gamma)^{2}+ (2-2\gamma)^2 }$. On the other
hand, by passing 
from $Z$ to $Z^{\mathrm{disc}}$, one may gain at most $2\gamma$ in
perimeter at $x = \pm\gamma$. Hence, we find 
$$P(Z^{\mathrm{disc}}, \Omega) \leq P(Z, \Omega) - 2 \sqrt{(2\gamma)^2+(2\!-\!2\gamma)^2} + 2\gamma\,.$$
This allows us to conclude that 
\begin{align}
P(Z^{\mathrm{disc}},\Omega)+\dissparam \calL^2(F^c {\setminus}Z^{\mathrm{disc}}) &\leq P(Z,\Omega) + \dissparam \calL^2( F^c {\setminus}Z) 
- 2\sqrt{(2\gamma)^2+(2-2\gamma)^2} + 2\gamma +4\dissparam
\gamma\nonumber\\
&<P(Z, \Omega) + \dissparam
\calL^{ 2}(F^c {\setminus}Z)\,, \label{contrad}
\end{align}
where the last strict inequality follows for $\gamma$ small
enough since
$$ - 2\sqrt{(2\gamma)^2+(2-2\gamma)^2} + 2\gamma +4 \dissparam
\gamma \to -4.$$ 
Before closing this discussion, let us point out that the forcing
$F$ from \eqref{fff} does not fulfill the
assumptions \eqref{hypF}. The argument above can however be reproduced
for a suitable smoothing of $F$ as well, at the expense of a somewhat
more involved notation.  
\paragraph{\bf Convexity.} In two dimensions, convexity is preserved 
by Stable Minimizing Movements. Again, to see this it is sufficient   to  check the preservation of convexity for the minimizers $Z$ of \eqref{c1}, 
 with $F^c$  convex. 
Assume by contradiction that a minimizer $Z$ is not convex and let 
$\overline{{\rm co}( Z)}$ be the closed convex hull of $ Z$. Owing
to \cite[Thm.\ 1]{Ferriero-Fusco}, we have that $ P(\overline{{\rm
    co}( Z)},\Omega )\leq
P(Z, \Omega)$. As $Z \subset F^c$ implies $\overline{{\rm co}(Z)}\subset \overline{{\rm co}(F^c)}
\equiv F^c $ and we conclude that 
$$P({\rm co}( Z),\Omega  )+ \dissparam \calL^{ 2}(F^c {\setminus} \overline{{\rm co}( Z)})\leq
P( Z, \Omega ) + \dissparam\calL^{2}(F^c {\setminus} Z).$$
\noindent
In particular, $\overline{{\rm co}( Z)}$ is a minimizer, too. This implies
that, if the  initial datum $Z_0$ for the whole evolutionary process  
is convex and the forcing term $F: [0,T] \rightrightarrows  \Omega$ is such
that $F^c(t)$ is convex for all $t\in (0,T]$, then any element 
$Z\in \SMM \infty$ is 
such that
$Z(t)$ is closed and convex for all $t\in[0,T]$.
\par
Note that we cannot apply the same argument in order to ensure that
{\it star-shapedness} with respect to some given point is preserved
along the evolution, for star-shaped rearrangements do not necessarily
decrease the perimeter, cf.\ \cite[Lemma 1.2]{Kawohl}.
\paragraph{\bf Symmetries.} If $Z_0$ and
$F^c(t)$ are balls (radially symmetric), then every  $Z\in \SMM \infty$ is a
ball for all $t\in (0,T]$ as well. 
This can be checked  by induction on 
 minimizers for problem  \eqref{c1}: Assume $Z^o$ and
$F^c$ to be radially symmetric. If $Z$ were not radially
symmetric one would strictly decrease the perimeter by redefining $Z$
to be
a ball with the same volume, included in $ F^c$. 
\par
\begin{figure}[h]
  \centering
  \pgfdeclareimage[width=95mm]{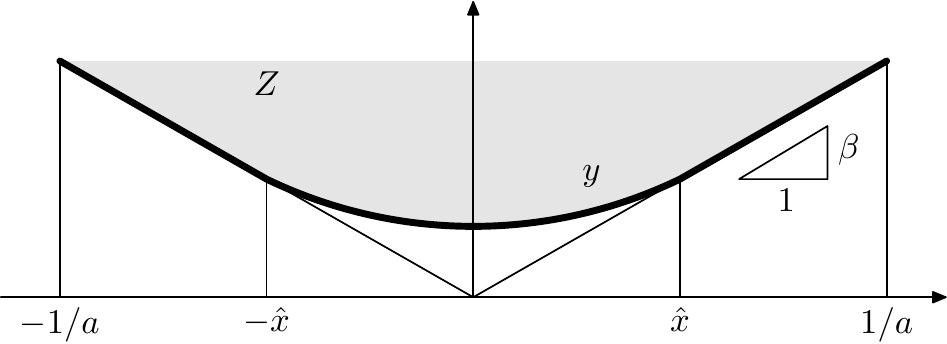}{regularity}
    \pgfuseimage{regularity}
  \caption{The $C^1$ competitor profile.}
\label{regularity}
\end{figure}
\begin{figure}[h]
  \centering
\pgfdeclareimage[width=145mm]{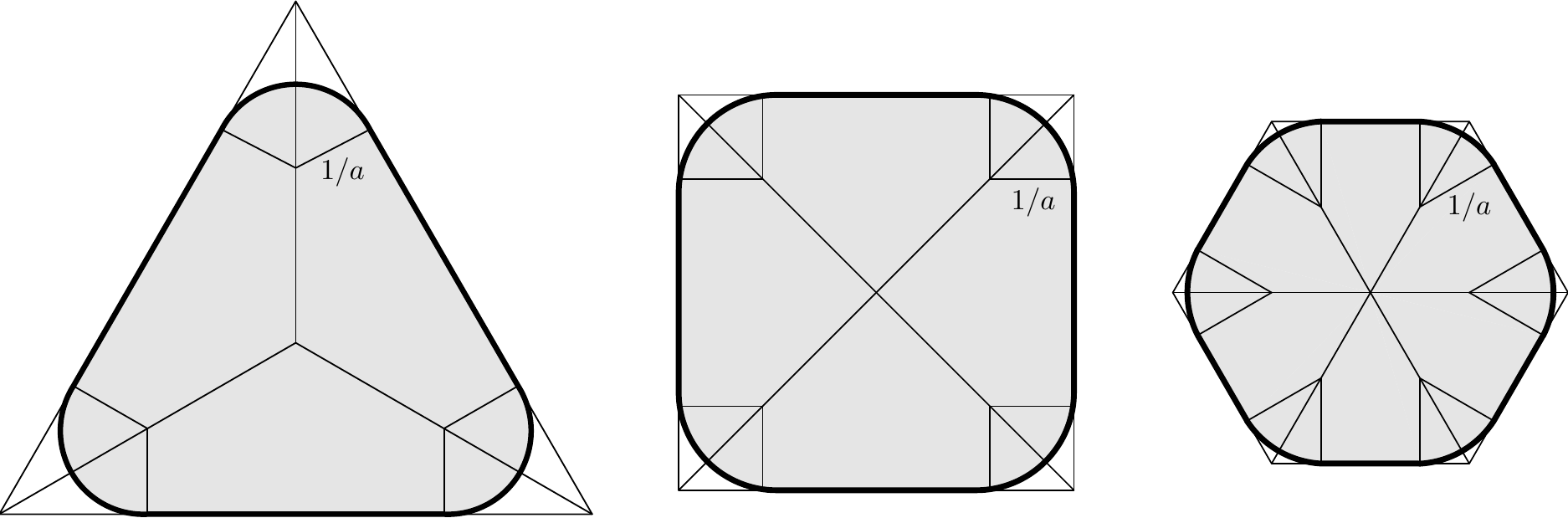}{hexagon}
    \pgfuseimage{hexagon}
  \caption{Solutions of the minimization problem \eqref{c1} with
    forcing $F^c$ being an equilateral triangle, a
    square, and a regular hexagon, respectively. }
\label{square}
\end{figure}

Analogously, other symmetries can be conserved along the
evolution.
For instance, let  $Z_0$  be
 symmetric with respect to a fixed hyperplane $\pi$ and suppose that the sets
 $F^c(t)$ have the same property for all $t\in (0,T]$.
  Then, any element 
$Z\in \SMM \infty$ 
is symmetric with respect to
$\pi$ as well.   We shall check this again at the  level of the time-incremental problem \eqref{c1}, 
supposing $Z^o$ and $F^c$ to be symmetric with respect to $\pi$.
If  $Z$ were not symmetric,  one could replace 
$Z$ with its Steiner symmetrization $Z^{\mathrm{s}}$  with respect to
$\pi$. This would  be admissible, since   $F^c$ is 
symmetric with respect to $\pi$. Moreover, one would have  that
$\calL^d(F^c{\setminus}Z^{\mathrm{s}})=\calL^d(F^c{\setminus}Z)$ and $P(Z^\mathrm{s}, \Omega)
\leq P(Z, \Omega)$. Figure \ref{square} shows some
examples of symmetric  minimizers. 
\paragraph{\bf Partial $\mathrm{C}^1$ regularity.} Regularity of the evolving set cannot be
expected in general, for it is easy to design forcing sets $F(t)$
resulting in reentrant corners  of the solution set  (i.e., points $x$ at the
boundary such that the set locally has  a cone of amplitude strictly larger than
$\pi$ and vertex at $x$), see
Figure \ref{f3} (right). On the other hand, in two dimensions, the set is smooth out of reentrant corners. 
 We will check this by considering the case of Cartesian
graphs. Let $ F^c$ be locally the epigraph of the piecewise affine
function $[0,1/\dissparam] \ni x \mapsto \beta|x|$ for $\beta>0$. We will check that the minimizer of \eqref{c1}
is a $\mathrm{C}^1$-set. In order to see this, we show that the minimizer $y \in W^{1,1}(0,1/\dissparam)$ 
of
\begin{equation}
\calF(y):=\int_{0}^{1/\dissparam} \left(\sqrt{1+(y'(x))^2}
    + \dissparam (y(x)-\beta x)\right)\,{\rm d} x 
\label{c2}
\end{equation}
under the conditions $y'(0+)=0$ and $ y(x) \geq
    \beta x $ for all $x \in (0,1/\dissparam)$
is indeed $\mathrm{C}^1$, see Figure \ref{regularity}. 
Minimum problem \eqref{c2} corresponds to problem \eqref{c1} under
the assumption that the minimizer $Z$ is symmetric w.r.t. the $y$
axis, under mild integrability assumptions. By assuming that the
optimal profile in \eqref{c2} is actually not in contact with the
constraint $\beta|x|$ in some (still unknown) interval $(-\hat x,\hat
x)$ for some $\hat x \in (0,1/a)$, one can consider variations of $\calF$
which are symmetric and compactly supported in $(-\hat x,\hat
x)$ in order to compute the Euler-Lagrange equation 
$$\frac{\mathrm{d}}{\mathrm{d}x}\frac{y'(x)}{ \sqrt{1+(y'(x))^2}} = \dissparam.$$ 
We can now solve for $y',$ taking into account
$y'(0+)=0$, and deduce that 
$$y'(x) = \frac{\dissparam x}{ \sqrt{1-\dissparam^2x^2}}\quad
\text{for all }x\in[0, \hat x)$$
which, by direct integration gives
$$y(x)=y(0)-\frac{1}{\dissparam}\sqrt{1-\dissparam^2x^2}.$$ 
This in particular entails that, independently of the opening
$\beta$, in case of no contact with the constraint $\beta|x|$ the
optimal profile is an arc of a circle with radius $1/a$. Note indeed
that the latter expression makes sense for $|x|\leq |\hat x| \leq 1/a$
only. 

In order to determine $y(0)$, we ask $\hat x$ to be a tangency
point between the optimal profile $y(x)$ and the constraint
$\beta|x|$. Indeed, such tangency must occur at some point in
$(0,1/a)$. If this were not the case, one could translate the profile
as $y(x)-k$ for $k>0$ up to tangency, which would contradict  
minimality. We can hence assume that $y(\hat x)=\beta\hat x$ and  $y'(\hat x)=\beta$
(otherwise this very argument could be repeated at $\hat x$, giving
rise to a contradiction). Moreover, $y(x) = \beta  x$ for all $x
\in (\hat x,1/a)$, for the only other options would be to have an arc
of radius $\dissparam$ in $(\hat x,1/\dissparam)$ as well, which would again contradict minimality. This gives 
$$\hat x = \frac{\beta}{\dissparam \sqrt{1+\beta^2}} \ \ \text{and} \ \
y(0)=\frac{1}{\dissparam}\sqrt{1+\beta^2}.$$
Note that $\hat x < 1/\dissparam$ for all $\beta>0.$ 
In particular, the candidate  optimal profile is 
\begin{equation}
y(x) =
-\frac{1}{\dissparam}\sqrt{1-\dissparam^2x^2}+\frac{1}{\dissparam}\sqrt{1+\beta^2}.\label{optimal}
\end{equation}
We now check that the $\mathrm{C}^1$ profile $y$ is optimal by comparing the value of  $\calF$ for $y$ with its value for the affine function
$\ell(x) = \beta x;$  the latter corresponds to a nonreentrant
corner  (i.e., a
point $x$ at the
boundary such that the complement of the set locally contains a cone of amplitude strictly larger than
$\pi$ and vertex at $x$).  Using that $y\equiv \ell$ 
on $(\hat x,1/\dissparam)$ we have
\begin{align*}
  \calF(y) -\calF(\ell)&= \int_0^{\hat x} \sqrt{1
  +\frac{\dissparam^2x^2}{1-\dissparam^2 x^2}}\,{\rm d}x  - \sqrt{1+\beta^2}\hat x + \dissparam \int_0^{\hat x} \left(
  -\frac{1}{\dissparam}\sqrt{1-\dissparam^2x^2}+\frac{1}{\dissparam}\sqrt{1+\beta^2} - \beta x
  \right) \, {\rm d} x\\
&= \int_0^{\hat x} \frac{\dissparam^2 x^2}{\sqrt{1-\dissparam^2 x^2}} \, {\rm d} x -
  \frac{\beta^3}{2\dissparam(1+\beta^2)} <\frac{\beta^3}{3\dissparam(1+\beta^2)^{3/2}} -
  \frac{\beta^3}{2\dissparam(1+\beta^2)} <0.
\end{align*}
This in particular shows that a nonreentrant corner at scale $1/a$
is not admissible and that the minimizer is $\mathrm{C}^1$ instead. On the
other hand, the above argument is  scale-invariant  and nonreentrant
corners are hence excluded at any scale. Note that the optimal profile
$y$ is not $\mathrm{C}^2$, since $y''(\hat x-)>0$. 

By combining this analysis with the  remark on preservation of
convexity, we can conclude that, in case $F^c$ is convex and piecewise $\mathrm{C}^1$, 
the minimizer of problem \eqref{c1} is globally $\mathrm{C}^1$, see also  Figures
\ref{f2} (left and right) and \ref{f3} (left). This remark  makes
the conclusions of Proposition \ref{prop:LDE} sharper, for in two space dimensions one
can choose $\mathfrak{a}=1/2$.

\paragraph{\bf Regular polygonal forcing.} 
In view of the above discussion, the minimizer $Z$ of  \eqref{c1} can be
explicitly determined in case $F^c$ is a regular
polygon and $1/\dissparam$ is smaller than half its side. Indeed, the preservation of symmetries implies that the
minimizer $Z$ shares the same symmetries of $F^c$, with rounded
corners of radius $1/\dissparam$, see Figure \ref{square}.

\begin{examples}[Nontrivial evolutions under compatibility condition \eqref{compatibility}] 
\label{ex:4.8}\upshape
Let $B_r$ denote the open ball  in $\R^d$  with radius $r\geq0$. By
imposing the   compatibility condition
\eqref{compatibility} one has that
$$ P(B_r,\Omega) = a  \calL(B_r) \quad \Leftrightarrow \quad
\omega_d r^{d-1} = a \frac{\omega_d}{d}r^d$$
where $\omega_d$ is the surface of the unit sphere in $\R^d$. Hence,
the only ball fulfilling \eqref{compatibility} has radius $d/a$. An
evolution under condition
\eqref{compatibility} and spherical symmetry is necessarily trivial: the ball of radius $d/a$
vanishes as soon as it is forced to evolve. Still, a first
nontrivial evolution example can be obtained by considering a disjoint
collection of balls of radius $d/a$.   For instance, the two-dimensional set 
\begin{equation}
Z(t) = \displaystyle\cup_{i=1}^{m(t)}B((i,0),2/a) \quad \text{for} \ \ m(t)
=M-\lfloor t \rfloor\label{zzz}
\end{equation}
for some $M\in \Nz$ (and $2/a<1/2$) (with $\lfloor\cdot\rfloor$ the Gauss-bracket), gives an
evolution corresponding to the forcing $F^c(t)=Z(t)$ and fulfills
\eqref{compatibility} for all times, see Figure \ref{palline}.
\begin{figure}[h]
  \centering
  \pgfdeclareimage[width=125mm]{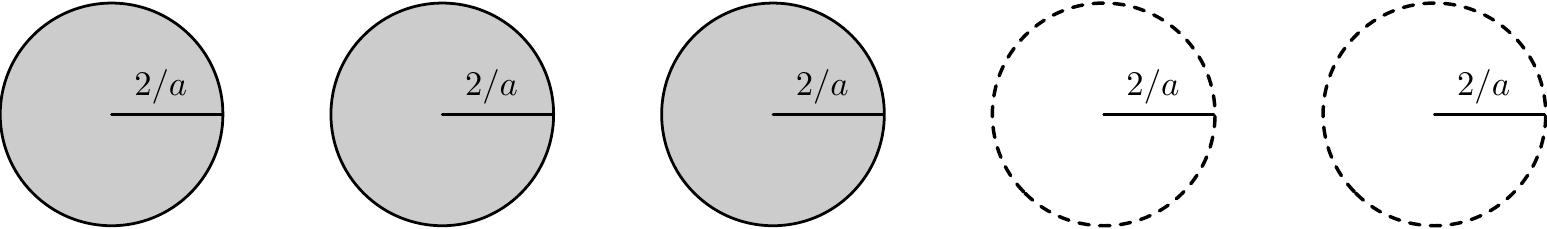}{palline}
    \pgfuseimage{palline}
  \caption{The evolution from \eqref{zzz} for $M=5$ and time $t=2$.}
\label{palline}
\end{figure}

In two dimensions, condition \eqref{compatibility} selects the unique minimizer
among each family of rounded polygonal shapes, see Figure
\ref{square}. By considering a collection of
disjoint smoothed polygons fulfilling condition \eqref{compatibility}
and balls of radius $2/a$ (hence fulfilling condition \eqref{compatibility}), one can again design a nontrivial evolution
in the spirit of Figure \ref{palline}.

Let us conclude by showing an example of a nontrivial evolution for a 
{\it connected} $Z(t)$ in two dimensions. Consider the smoothed square in the
middle of Figure \ref{square}. Elementary algebra shows that the only
smoothed square fulfilling condition \eqref{compatibility} is
inscribed in the square of side
$(2+\sqrt{\pi})/a$. In particular, the flat portion of each side
measures $\sqrt{\pi}/a\sim 1.77/a$. Let us now consider the union of
discs of radius $1/a$ and the smoothed square, as in Figure
\ref{mickey}. This union can be realized by still fulfilling condition
\eqref{compatibility}, as long as the positioning of each extra disc is such that
the gain in perimeter equates $a$-times the gain in
area. By calling $\alpha$ the angle at the center of the disk which
identifies the arc cut by the side of the smoothed square, the
aforementioned equality reduces to
$$ 2\pi - 
\alpha - 4 \sin(\alpha/2)- \sin(\alpha)=0$$
(note that it is independent of $a$), whose unique solution in
$[0,\pi]$ is $\tilde \alpha\sim 2.005$. Note that the cord of the disk
of radius $1/a$
corresponding to $\tilde \alpha$ has length $2\sin(\tilde \alpha/2)/a
\sim 1.687/a$, which is strictly shorter than the flat portion of each
side of the smoothed square. 
\begin{figure}[h]
  \centering
  \pgfdeclareimage[width=155mm]{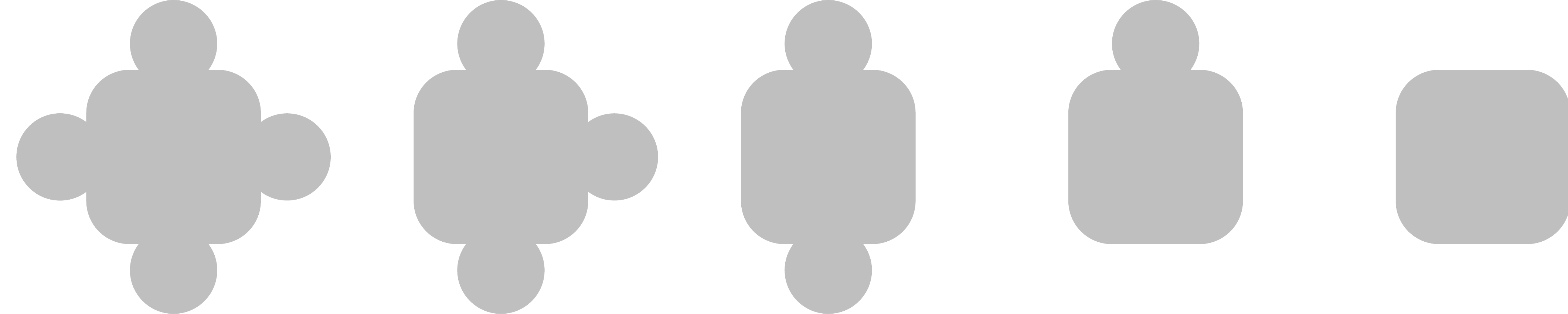}{mickey}
    \pgfuseimage{mickey}
  \caption{An evolution of connected sets fulfilling the
    compatibility condition
    \eqref{compatibility}, time flows from left to right.}
\label{mickey}
\end{figure}
All configurations in Figure \ref{mickey} hence fulfill 
condition \eqref{compatibility}. Moreover, they are stable (with respect
to the forcing corresponding to the interior of their complement), for
they fulfill the interior ball condition with balls of radius
$1/a$. An evolution as depicted in Figure  \eqref{mickey} (left to right) can hence be
realized by suitably prescribing the forcing. At all times, such
evolution fulfills  \eqref{compatibility}.
\end{examples}.
%
\subsection{A numerical test}
%
In order to illustrate the above discussion, we provide some numerical
evidence in a planar setting. We assume $\Omega =
(-3,3)^2$ and consider  an initial state $Z^o\subset \Omega$ such that
$F^c \subset Z^o$ where 
$$ F=\{(x,y) \in [0,1]\times [-1,1]  \ : \
 |y|  > v(x) \} \cup \{(x,y) \in \R^2 \ : \  x<0 \
\text{or} \ x>1 \},$$ where $v:[0,1] \to (0,1]$ is a given function,
different for each  numerical example.  In the following,
we seek  minimizers $Z$ of  the
incremental problem \eqref{c1} of the form 
$$
Z=\{(x,y) \in 
[0,1]\times [-1,1]  \ : \
|y| \leq u(x) \}$$
 for some optimal profile $u:[0,1]\to [0,1]$ to be determined 
 such that $u(x)\leq v(x)$ for all $x\in[-1,1]$, 
 which is in accordance with the brittle constraint $Z \subset F^c$.  
 Note that, given the discussion on symmetry
 from Subsection \ref{ss:ulisse},  assuming $Z$ to be symmetric with respect
 to $\{y=0\}$ is not restrictive, since $F^c$ also is. Moreover, owing to the discussion leading to \eqref{optimal},
 it is not restrictive to assume that $[0,1] \times \{0\} \subset Z$
 (that is, $Z$ can actually be described by the profile $u$) as long as
 the optimal profile $u$ fulfills $$u(x) \geq \sqrt{(a^2 - (x{-}a)^2)^+}+
 \sqrt{(a^2-(x{-}1{+}a)^2)^+}\quad \forall x \in [0,1]$$ (that is, if $Z$ contains the two
 balls of radius $a$ centered in $(a,0)$ and $(1{-}a,0)$), which happens to be the case for all
 computations below.

The problem is discretized  in space by partitioning the  domain $[0,1]$ of the variable $x$ as $0=x_0<x_1<\dots<x_N=1$ with
 $x_i=i/N$ and $N=100$  
 and by approximating $v$ via its piecewise
affine interpolant on the partition, taking the values
$v(i/N)=:v_i$. In particular, we look for  $u$ piecewise affine with
$u(i/N)=:u_i$ minimizing
\begin{align*}
(u_0,\dots,u_N) &\mapsto P(Z,\Omega) + \dissparam {\mathcal L}^2(F^c{\setminus} Z) \\
&
= 2\left(u_0+u_N + \sum_{i=1}^N\sqrt{(u_i{-}u_{i-1})^2+N^{-2}} +
\dissparam \sum_{i=1}^N N^{-1} |u_i- v_i| \right) \\
&\text{
under the constraints $0\leq u_i\leq v_i$ for $i=1,\dots,N$. }
\end{align*}
This is a
strictly convex minimization problem under convex constraints. In the
following, we solve it by
using the {\tt fmincon} tool of Matlab for different
choices of the function $v$. In all figures, we depict the {\it portions} of
the minimizer
$Z$ (light color) and of the
forcing set $F$ in $[0,1]^2$ (dark color). The reader should
however keep in mind that both forcing and minimizer are actually symmetric along $\{y=0\}$. 

\begin{figure}[h]
  \centering
  \pgfdeclareimage[width=75mm]{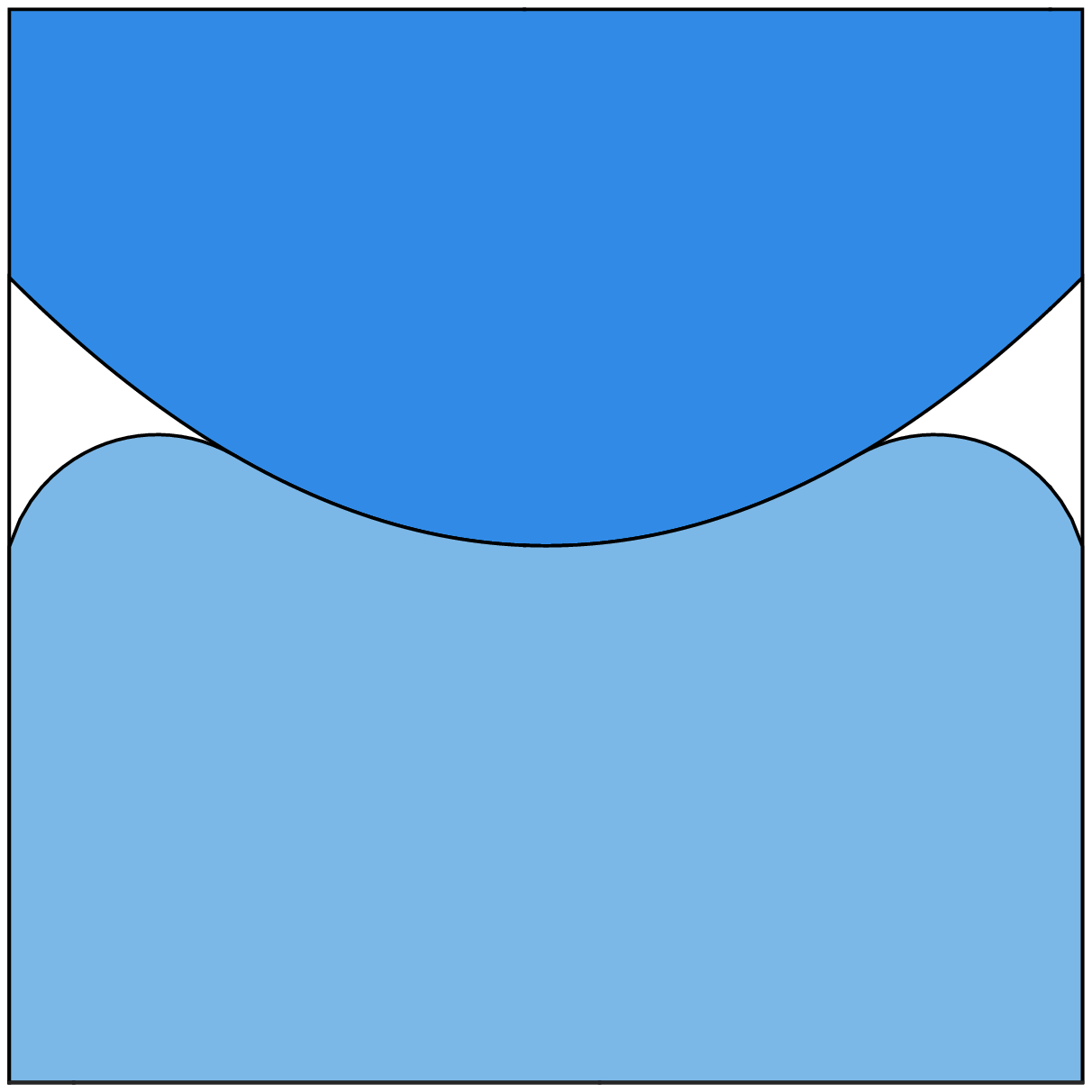}{f11}
  \pgfdeclareimage[width=75mm]{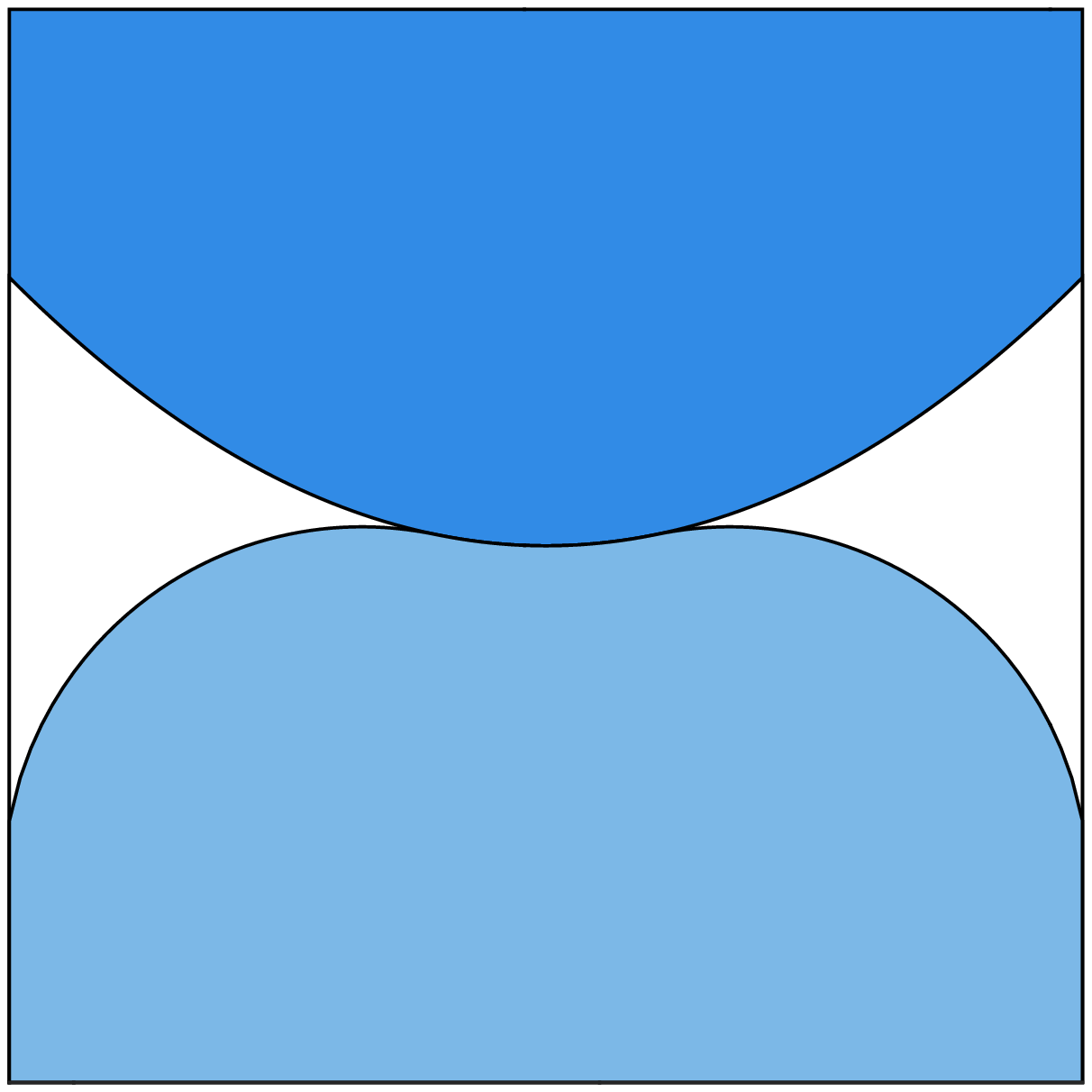}{f12}
  \pgfuseimage{f11}\hspace{-9mm}
  \pgfuseimage{f12}
  \caption{The effect of changing the parameter $\dissparam$. The two solutions
    correspond to  $v(x) = (x{-}1/2)^2+1/2$ for $\dissparam=7$ (left) and
    $\dissparam=3$ (right). The top adhesion zone is
    smaller for smaller $\dissparam$. Note that the parts of the
    boundary of $Z$ which are not in contact with $F^c$ are arcs of
    circles with radius $1/a$ (recall that $\dissparam$ is different in
    the two figures), as predicted in Subsection
    \ref{ss:ulisse}.}
  \label{f1}
\end{figure}

Figure \ref{f1} corresponds to the choice  $v(x) = (x{-}1/2)^2+1/2$ and 
illustrates the effect of changing the parameter
$\dissparam$. A smaller value of $\dissparam$ favors a shorter perimeter at the
expense of a larger distance from $F^c$. Correspondingly, the top adhesion zone, namely the
points where $u\equiv v$, is smaller for smaller $\dissparam$.

Let us mention that, in case $\dissparam<2$ one can prove that
$u_0=u_N=0$ which, as mentioned above, may well be not admissible. In order to avoid this pathology,
in all the following simulations, the parameter $\dissparam$ will be always
chosen to be $5$. Correspondingly, in all simulations the optimal
profile $u$ is everywhere well separated from the $y=0$ axis. 

\begin{figure}[h]
  \centering
  \pgfdeclareimage[width=75mm]{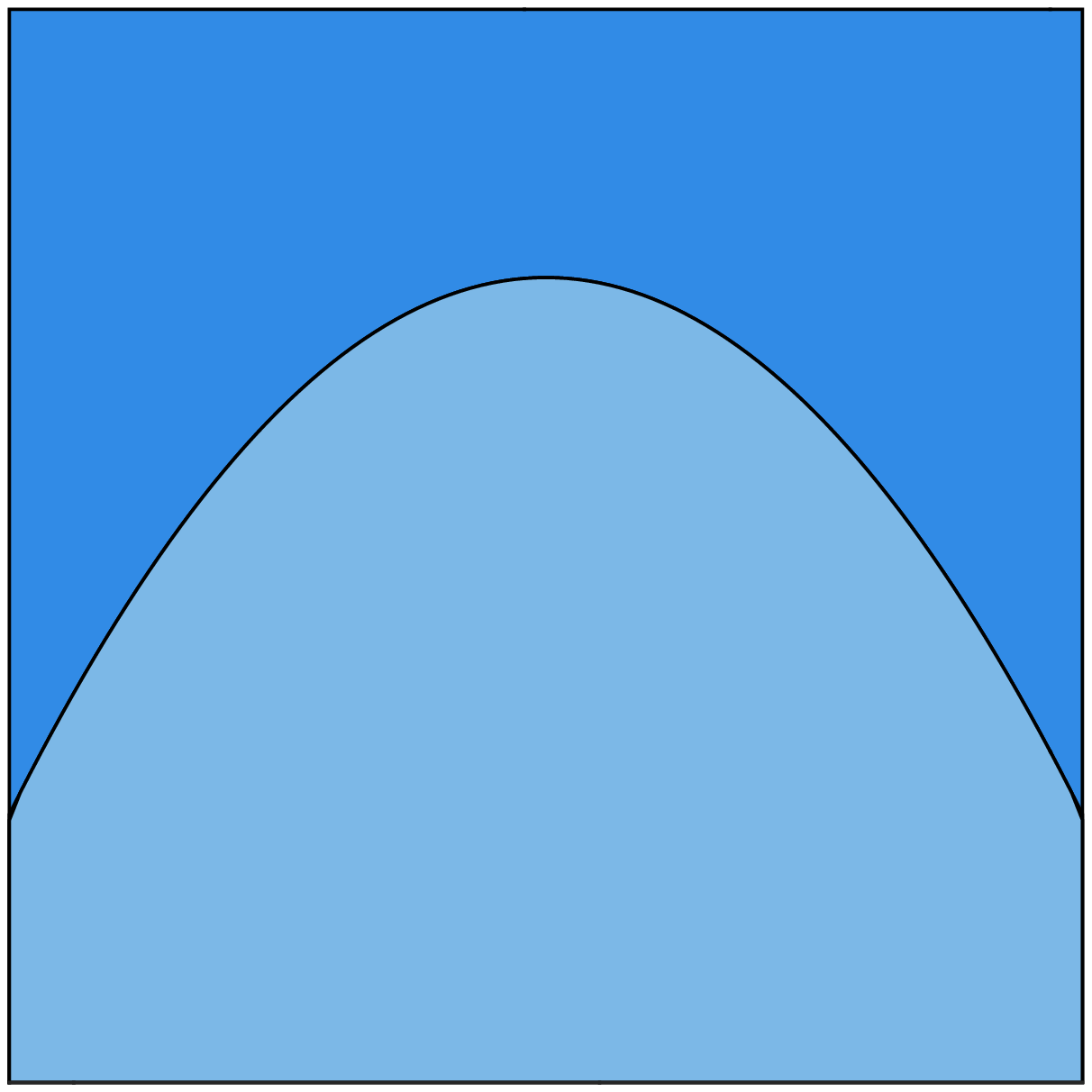}{f21}
  \pgfdeclareimage[width=75mm]{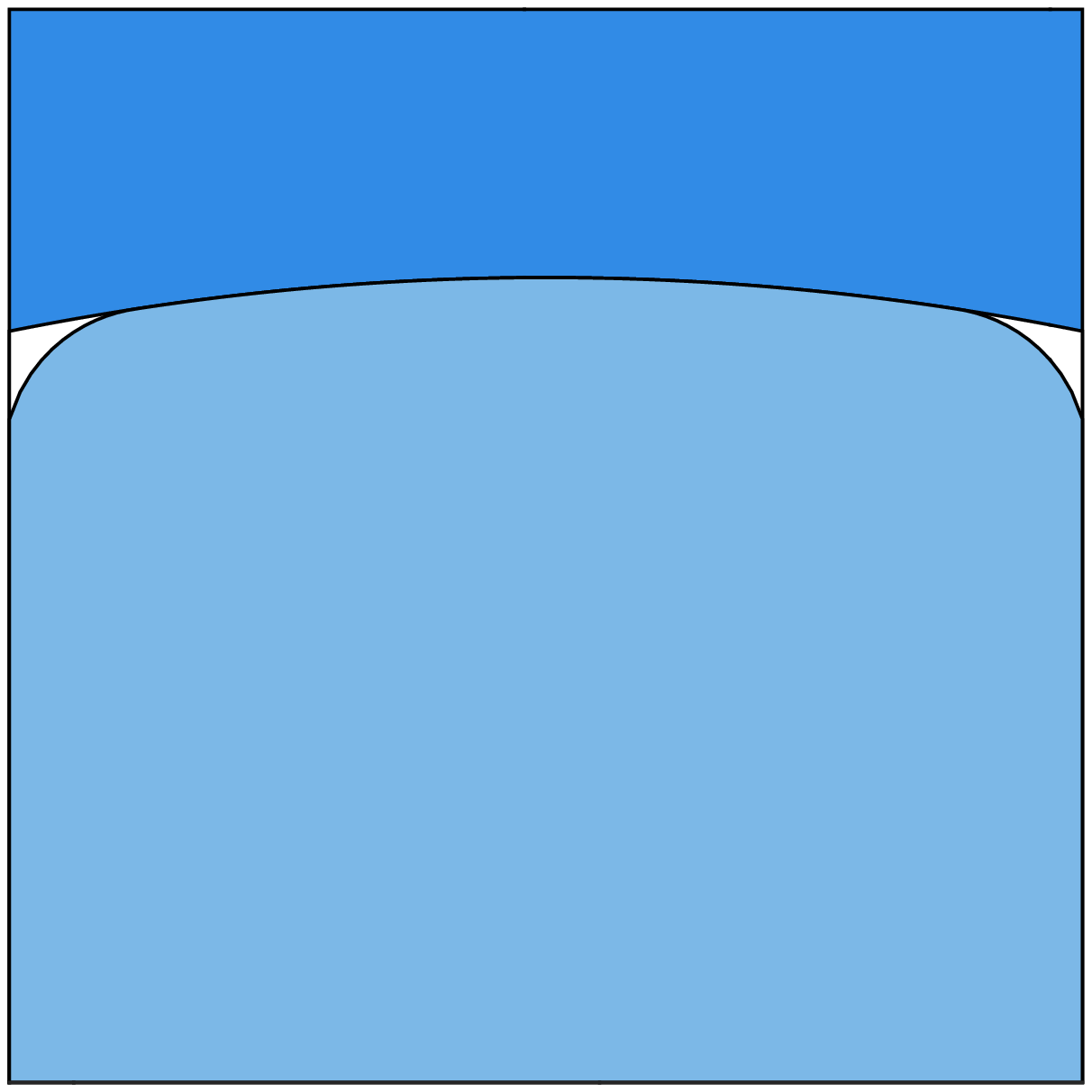}{f22}
  \pgfuseimage{f21}\hspace{-9mm}
  \pgfuseimage{f22}
  \caption{Convex forcing  $F^{c}$. The two solutions
    correspond to  $v(x) = 3/4 - \beta(x{-}1/2)^2$ for $\beta=2$ (left) and
    $\beta=1/5$ (right). The minimal set $Z$ is convex. }
  \label{f2}
\end{figure}

Figure \ref{f2} follows by letting  $v(x) = 3/4 - \beta(x{-}1/2)^2$  
along with two different choices of the parameter $\beta$ and is meant to
illustrate the convexity of the evolution in presence of a convex
forcing  $F^{c}$,  see Subsection 4.1. In particular, the minimizer $Z$ is
convex. One observes that the optimal profile $u$ detaches from $v$
even in the convex case. This is indeed the case also in Figure
\ref{f2} left, where nonetheless the detachement is not visible due to
the scale. 

\begin{figure}[h]
  \centering
  \pgfdeclareimage[width=75mm]{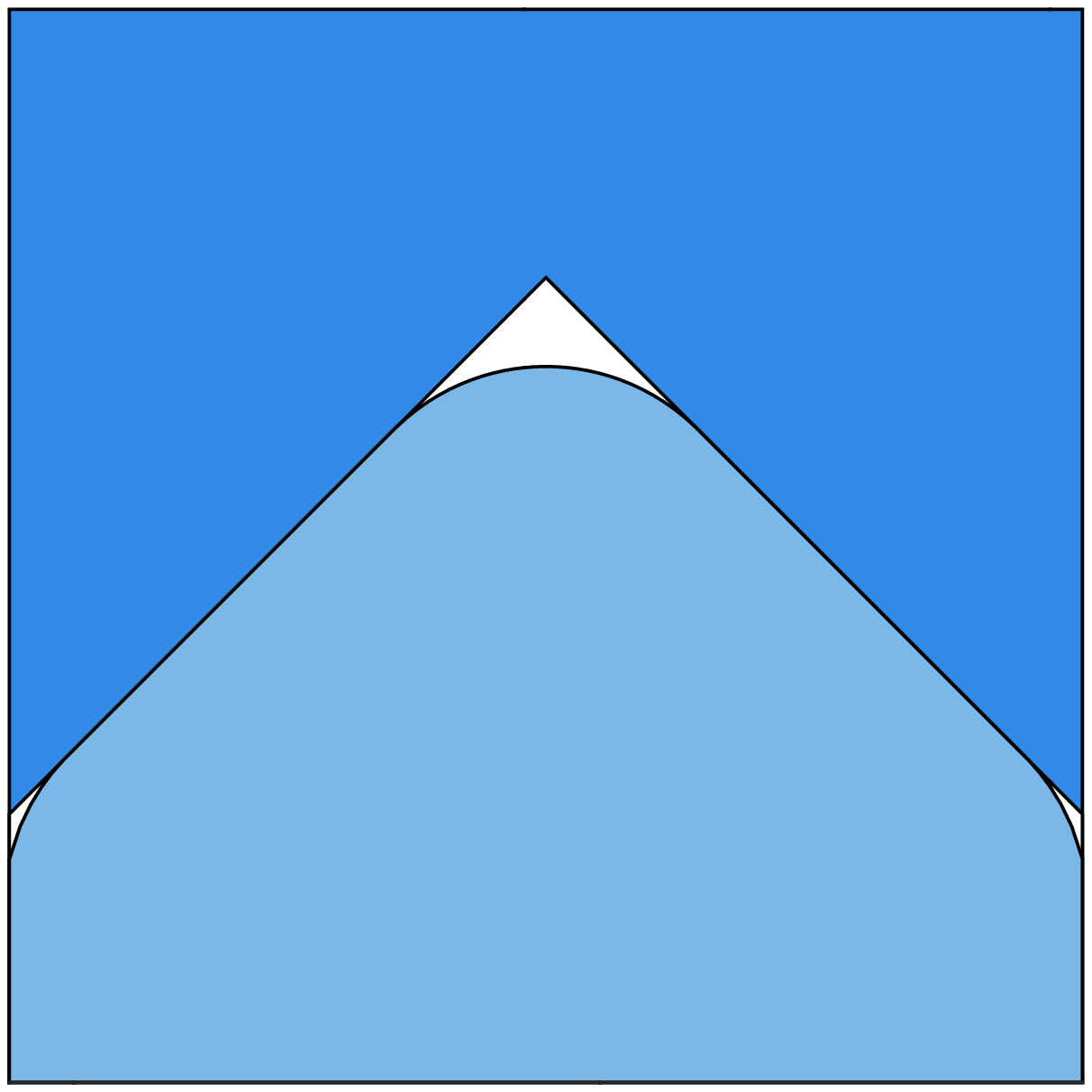}{f31}
  \pgfdeclareimage[width=75mm]{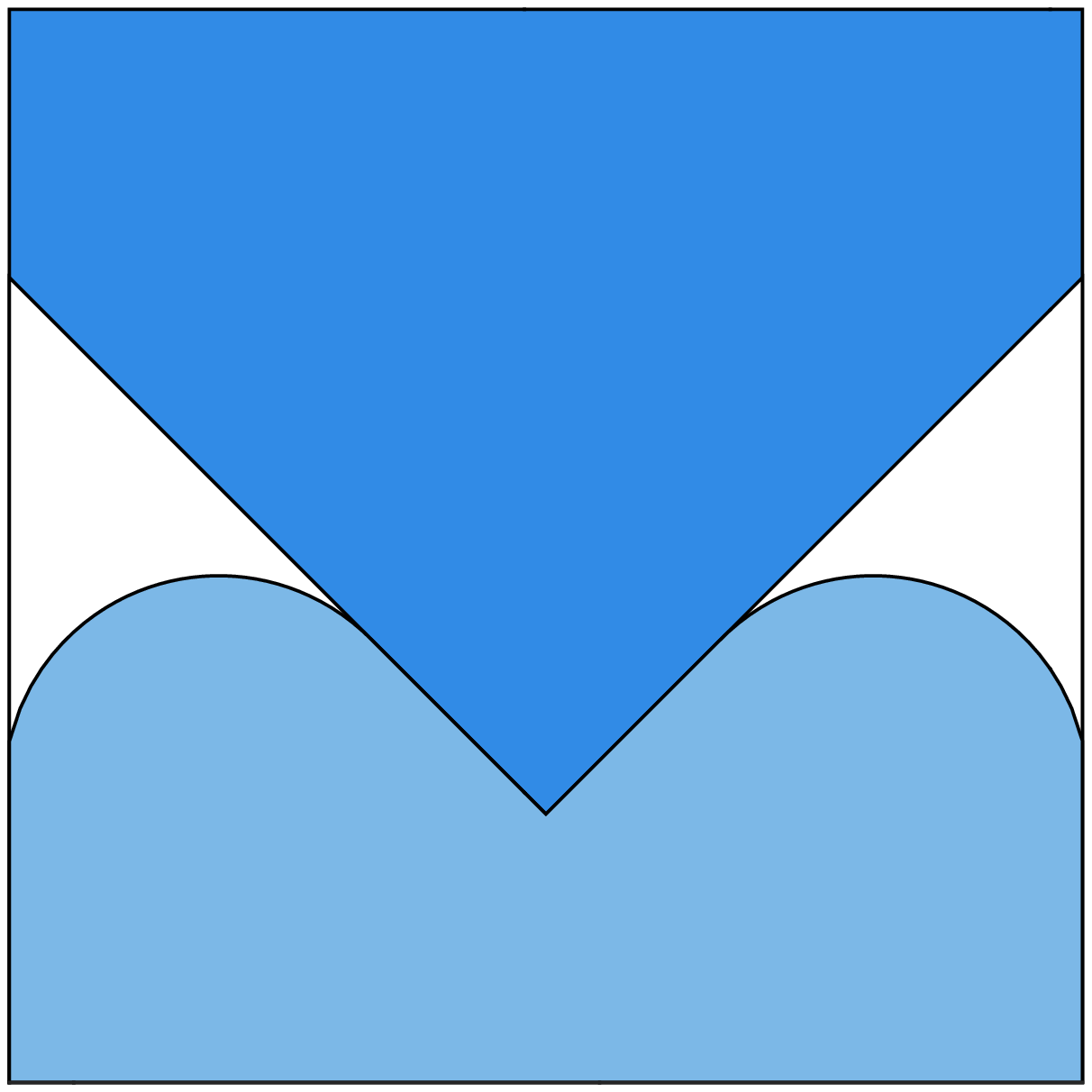}{f32}
  \pgfuseimage{f31}\hspace{-9mm}
  \pgfuseimage{f32}
  \caption{Partial $C^1$ regularity. The two solutions
    correspond to  $v(x) = 3/4 - |x{-}1/2|$ (left) and  $v(x) = 1/4
    + |x{-}1/2|$ (right).}
  \label{f3}
\end{figure}

Figure \ref{f3} corresponds to the choices 
 $v(x) = 1/2 \pm (|x{-}1/2|-1/4)$ and illustrates
the partial regularity of the solution. Note that nonsmooth 
boundary points occur in connection with nonconvex forcings  $F^{c}$.  

\begin{figure}[h]
  \centering
   \pgfdeclareimage[width=75mm]{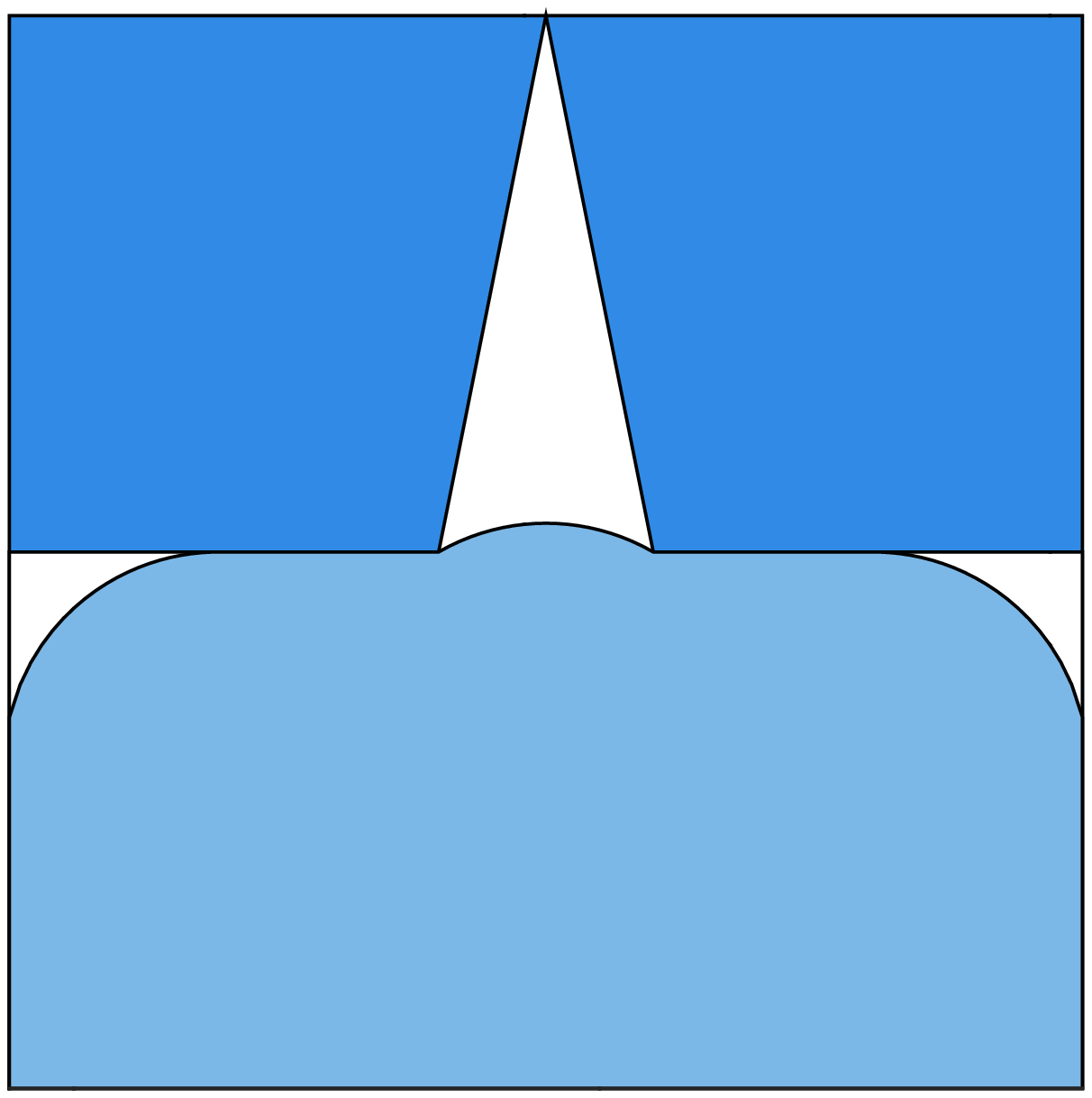}{f41}
  \pgfdeclareimage[width=75mm]{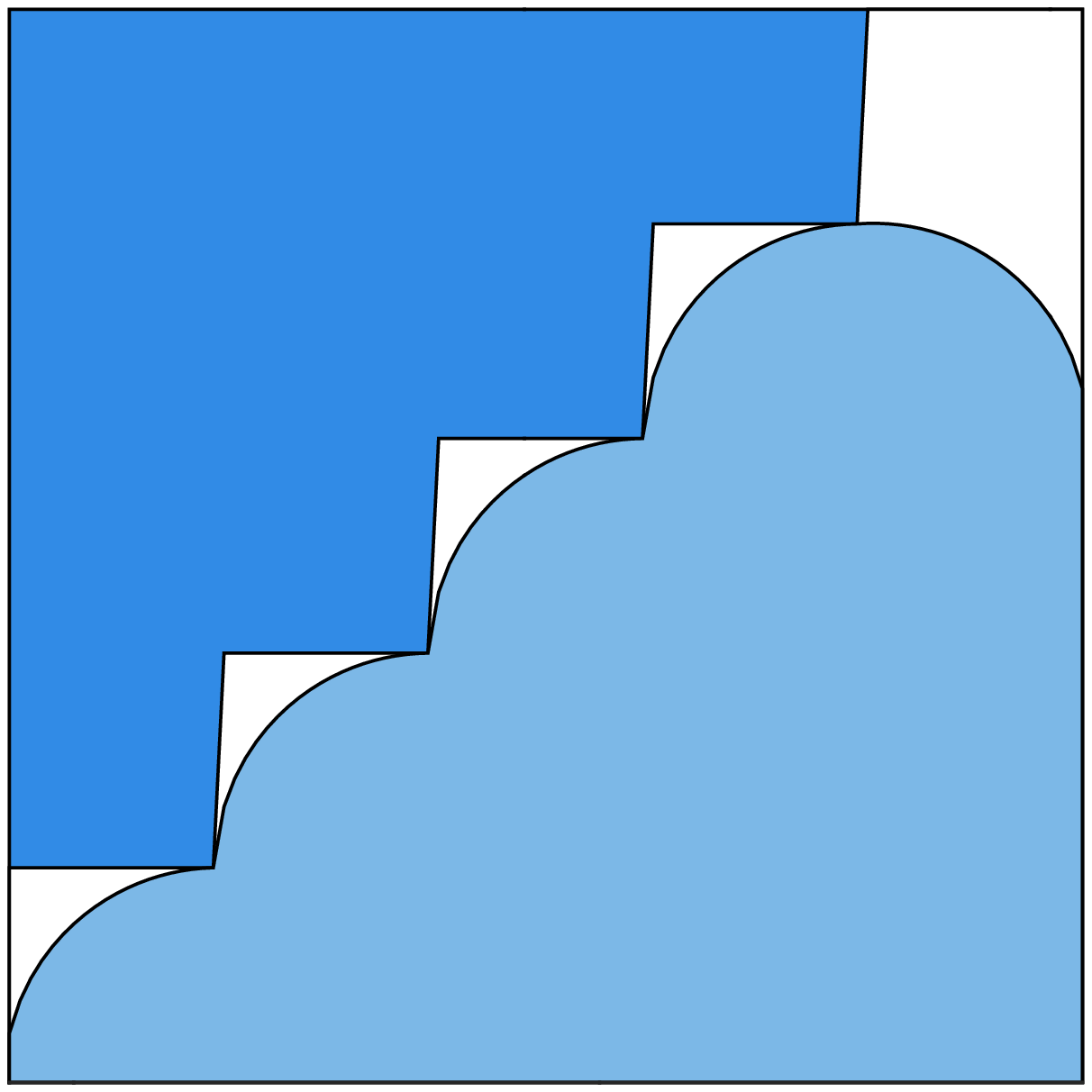}{f42}
  \pgfuseimage{f41}\hspace{-9mm}
  \pgfuseimage{f42}
  \caption{Extreme configurations. The solution
for  $v(x)=\max\{1-5|x{-}1/2|,1/2\}$ (left) and  $v(x)= \lfloor 5x\rfloor/5+1/5$ (right).}
  \label{f4}
\end{figure}

Figure \ref{f4} illustrates some special situations. On the left, the solution
for  $v(x)=\max\{1-5|x{-}1/2|,1/2\}$. In this case, the optimal
profile is such that $u(1/2)>3/4$. Note that the same holds for any
opening of the cone in  $F^c$, whatever small.
On the right, the solution for  $v(x)= \lfloor 5x\rfloor/5+1/5$. The profiles $u=v$ touch at the points
  $x=1/5, \, 2/5, \, 3/5$, and $4/5$ only.  Note that all
  nonstraight portions of the boundary of $Z$ are arcs of radius
  $1/\dissparam$. 

%
\section{Proofs of Theorems \ref{thm:4.1}, \ref{thm:4.2}, and \ref{thm:4.3}}
\label{s:6}
%
\noindent
We start by the \underline{\bf Proof of Theorem \ref{thm:4.1}:} For  $k\in \N\cup \{\infty\}$ fixed, let $(\tau_j)_j$ be a vanishing sequence 
and let $(\pwc z{\tau_j})_j$ be the characteristic functions of the sets $(\pwc Z{\tau_j})_j$. Since the functions 
$\pwc z{\tau_j}$ are nonincreasing, it is immediate to check that \[
\exists\, C>0 \ \forall\, j \in \N \, : \qquad \| \pwc z{\tau_j}\|_{\mathrm{BV}([0,T];L^1(\Omega))} \leq C\,.
\]
Furthermore, 
from \eqref{energy-bound-interp} we get that 
\[
\exists\, C>0 \ \forall\, j \in \N \, : \qquad \| \pwc z{\tau_j}\|_{L^\infty(0,T;\mathrm{SBV}(\Omega;\{0,1\}))} \leq C\,.
\]
We are now in the position to apply a Helly-type compactness result (cf.\ e.g.\ \cite[Thm.\ 3.2]{MaiMie05EREM}), and conclude that 
\[
\begin{aligned}
&
\exists\, z \in L^\infty(0,T;\mathrm{SBV}(\Omega;\{0,1\})) \cap \mathrm{BV}([0,T];L^1(\Omega)),
\\
& 
 \quad \text{with } 
z(\cdot, x) \text{ nonincreasing on } [0,T] \ \foraa\, x \in \Omega, \text{ such that }
\\
&
 \pwc z{\tau_j}(t) \weaksto z(t) \quad \text{ in } \mathrm{SBV}(\Omega;\{0,1\}) \quad \text{for every } t \in [0,T].
 \end{aligned}
\]
Then, the curve $Z: [0,T]\to \Spx$ defined by
$Z(0)=Z_0$ and 
  $Z(t): =\{ x\in \Omega\, : \ z(t)=1\}$ for all $t\in (0,T]$  is in $\GMM k$. 
  \par
  In order to prove convergence \eqref{refined-convs-1} for the  piecewise constant, right-continuous  interpolants $(\upwc z{\tau_j})_j$ we may argue in this way: by the above Helly argument, 
  there exists
$\underline z \in L^\infty(0,T;\mathrm{SBV}(\Omega;\{0,1\})) \cap \mathrm{BV}([0,T];L^1(\Omega))$ such that, 
   up to a further subsequence,   $ \upwc z{\tau_j}(t) \weaksto \underline{z}(t) $ in $ \mathrm{SBV}(\Omega;\{0,1\}) $ for all $t\in [0,T]$. Let $J$ be the  union of the jump sets of $z$ and $\underline z$: arguing for instance as in the proof of \cite[Thm.\ 4.1]{ThoRouPan},  
    it can be checked that $z(t) = \underline{z}(t)$ for all $t\in [0,T]\setminus J$. 
   Convergences \eqref{refined-convs-2} and  \eqref{refined-convs-3}
ensue from standard weak and strong compactness arguments. 
   \par
  Let us now pick $Z\in \GMM k $ approximated by a sequence $(\pwc Z{\tau_j})_j$  in the sense of 
\eqref{charact}. In order to  show that $Z$ satisfies the stability condition \eqref{stab}, we will pass to the limit in its discrete 
version \eqref{discr-stab-interp}, satisfied by the functions  $\pwc Z{\tau_j}$, by verifying the so-called \emph{mutual recovery sequence} condition from \cite {MRS06}. Namely, for every  fixed $t\in (0,T]$ and every admissible competitor $\widetilde Z \in \Spx$  for \eqref{stab}, with associated characteristic function $\tilde z$, we will exhibit a sequence $(\tilde z_j)_j\subset  \SBV(\Omega;\{0,1\})$ such that 
\begin{equation}
\label{MRS}
\limsup_{j\to\infty} \left( \calE_k(t,\tilde{z}_j) {-}\calE_k(t,\pwc z{\tau_j}(t)) {+} \calD(\pwc z{\tau_j}(t),\tilde{z}_j) \right) \leq  \calE_k(t,\tilde{z}) {-}\calE_k(t,z(t)) {+} \calD( z(t),\tilde{z}) \,. 
\end{equation}
The construction of the sequence $(\tilde z_j)_j$ is slightly adapted from the proof of 
\cite[Prop.\ 5.9]{RosTho12ABDM},  which in turn follows the steps of 
\cite[Lemma 2.13]{Thom11QEBV}. Hence, in the following  lines we shall refer to   \cite{RosTho12ABDM} for some details. 
First of all, we suppose that 
$ \calD( z(t),\tilde{z}) <\infty$, whence $\tilde z \leq z(t) $ a.e.\ in $\Omega$,   and, if  $k=\infty$, that $\calJ_\infty(t, \tilde z)=0$, so that  
$ \calE_\infty(t,\tilde{z}) <\infty$; otherwise, there is nothing to prove. 
Along the foosteps of \cite{RosTho12ABDM}, we set
\begin{equation}
\label{recovery-j}
\tilde{z}_j: = \tilde z \chi_{A_j} +\pwc z{\tau_j}(t) (1{-} \chi_{A_j}),  \quad \text{where } A_j: = \{ x \in \Omega\, : \ 0\leq \tilde z(x) \leq \pwc z{\tau_j}(t,x) \} \,.
\end{equation}
This way, we ensure that 
\begin{equation}
\label{1st-props-tildez-k}
\tilde{z}_j \in \{0,1\} \quad\text{and}\quad 0 \leq \tilde{z}_j \leq \pwc z{\tau_j}(t)  \qquad \aein\ \Omega\,.
\end{equation}
Furthermore, arguing in the very same way as in the proof of \cite[Prop.\ 5.9]{RosTho12ABDM}, 
where \cite[Thm.\ 3.84]{AmFuPa05FBVF}   on the decomposition of $\mathrm{BV}$-functions  is applied, we can show that 
$\tilde z_j \in \mathrm{SBV}(\Omega;\{0,1\})$. We now split the proof of 
\eqref{MRS} in $3$ steps:
\begin{enumerate}
\item Since 
\begin{equation}
\label{Lq-conv}
\pwc z{\tau_j}(t)\to z(t) \qquad \text{in $L^q(\Omega)$ for all $1\leq q <\infty$},
\end{equation}
 and 
$\tilde z \leq z(t) $ a.e.\ in $\Omega$, it is immediate to infer that 
$\tilde{z}_j \to \tilde z$ a.e.\ in $\Omega$ as $j\to\infty$, which improves to
\begin{equation}
\label{strong-converg}
\tilde{z}_j \to \tilde z \quad \text{in } L^q(\Omega) \quad \text{for all } 1\leq q <\infty,
\end{equation}
as $\tilde{z}_j \in \{0,1\}$ a.e.\ in $\Omega$. Combining \eqref{1st-props-tildez-k} with
\eqref{Lq-conv} and  
\eqref{strong-converg} we ultimately conclude that 
\begin{equation}
\label{conv-dissip}
\lim_{j\to\infty} \calD(\pwc z{\tau_j}(t),\tilde{z}_j) = \calD( z(t),\tilde{z}) \,.
\end{equation}
\item For $k\in \N$, we observe that 
\begin{subequations}
\label{limsupJ}
\begin{equation}
\label{limsupJk}
\begin{aligned}
\limsup_{j\to\infty} \left( \calJ_k(t,\tilde{z}_j) {-}\calJ_k(t,\pwc z{\tau_j}(t))\right) = 
\limsup_{j\to\infty} \int_{\Omega} k f(t) (\tilde{z}_j {-}  \pwc z{\tau_j}(t)) \ddx x 
 & =  \int_{\Omega} k f(t) (\tilde{z}{-}  z(t)) \ddx x \\ &  = \calJ_k(t,\tilde{z}) {-}\calJ_k(t,z(t))
 \end{aligned}
\end{equation}
thanks to \eqref{Lq-conv} and  
\eqref{strong-converg}.
For $k=\infty$,
we first of all observe that, since $\sup_{j\in \N} \calE(t,\pwc z{\tau_j}(t)) \leq C $ by \eqref{charact}, 
 we have $\calJ_\infty (t,\pwc z{\tau_j}(t)) =0$, i.e.\ $f(t,\cdot) \pwc z{\tau_j}(t,\cdot)=0$ a.e.\ in $\Omega$, for every $j\in \N$. Then, $0\leq  \calJ_\infty(t,z(t)) \leq \liminf_{j\to\infty} \calJ_\infty (t,\pwc z{\tau_j}(t)) =0$. Furthermore, 
 since 
$ \tilde{z}_j \leq \pwc z{\tau_j}(t) $ a.e.\ in $\Omega$ thanks to \eqref{1st-props-tildez-k}, we ultimately conclude that $\calJ_\infty (t,\tilde{z}_j) =0$ for all $j\in \N$. Therefore,
\begin{equation}
\label{limsupJ-infty}
\limsup_{j\to\infty} \left( \calJ_\infty(t,\tilde{z}_j) {-}\calJ_\infty(t,\pwc z{\tau_j}(t))\right) =
0=  \calJ_\infty(t,\tilde{z})  -  \calJ_\infty(t,z(t))
\end{equation}
\end{subequations}
(recall that we have supposed right from the start that $\calJ_\infty(t,\tilde{z})=0$). 
\item Finally, it can be shown that 
\begin{equation}
\label{limsup-perims}
\limsup_{j\to\infty} \left( |\mathrm{D} \tilde{z}_j|(\Omega) {-}  |\mathrm{D} \pwc z{\tau_j}(t)|(\Omega)\right) \leq |\mathrm{D} \tilde{z}|(\Omega) {-}  |\mathrm{D} z(t)|(\Omega)
\end{equation}
by repeating the very same arguments as in the proof of \cite[Prop.\ 5.9]{RosTho12ABDM}. 
\end{enumerate}
Combining \eqref{conv-dissip}, \eqref{limsupJ}, and \eqref{limsup-perims}, we infer \eqref{MRS} for $k\in \N$ and $k=\infty$, 
which concludes the proof of the stability condition \eqref{stab}, and thus of Thm.\ \ref{thm:4.1}. 
\QED
\paragraph{\bf \underline{Proof of Theorem \ref{thm:4.2}:}} Let $k\in \N$ be fixed, 
let $Z\in \SMM k$, and let $(\pwc z{\tau_j})_j$ converge to $z=\chi_Z$ as in \eqref{charact}. 
We can pass to the limit in the upper energy-dissipation estimate in \eqref{enbd} by observing that
\[
\begin{aligned}
&
\liminf_{j\to\infty} \left(\calE_k(\pwc{\sft}{\tau_j}(t), \pwc Z{\tau_j}(t)) {-}\calE_k(t,z(t))\right)   \\ & \geq \lim_{j\to\infty}  \left(\calE_k(\pwc{\sft}{\tau_j}(t), \pwc Z{\tau_j}(t)) {-} \calE_k(t, \pwc Z{\tau_j}(t))\right) + \liminf_{j\to\infty} \left(\calE_k(t, \pwc Z{\tau_j}(t)) {-}\calE_k(t,z(t))\right)  =: l_1+l_2 \geq 0,
\end{aligned}
\]
where we have used that 
\[
l_1=\lim_{j\to\infty}  \int_\Omega k \left( f(\pwc{\sft}{\tau_j}(t),x) {-} f(t,x)\right) 
 \pwc z{\tau_j}(t,x)) 
 \ddx x \leq k \lim_{j\to\infty}  \int_{t}^{\pwc{\sft}{\tau_j}(t)} \|\partial_t f(s) \|_{L^1(\Omega)} \ddx s  =0,
 \]
 due to \eqref{assFsmooth}, 
while, by convergence \eqref{charact} we have 
\[
l_2 =  \liminf_{j\to\infty}\left( |\mathrm{D}\pwc z{\tau_j}(t)|(\Omega) {-} |\mathrm{D} z(t)|(\Omega) \right) \geq 0\,.
\]  
Thanks to \eqref{charact}, we also  have 
$ \liminf_{j\to\infty} \calD(Z_0,\pwc Z{\tau_j}(t)) \geq  \calD(Z_0, Z(t))  $, and  \eqref{refined-convs} gives 
\[
\begin{aligned} 
\lim_{j \to\infty}
\int_{0}^{\pwc \sft{\tau_j}(t)}\partial_t\calE_k(r,\upwc Z{\tau_{j}}(r))\,\mathrm{d}r 
 & = \lim_{j \to\infty}
\int_{0}^{\pwc \sft{\tau_j}(t)}  \int_\Omega k \partial_t f(r,x)  \upwc z{\tau_{j}}(r,x)  \ddx x \,\mathrm{d}r 
\\
 & = \int_0^t \int_\Omega  k\partial_t f(r,x) z(r,x)  \ddx x \ddx r   =\int_0^t \partial_t \calE_k(r,Z(r)) \ddx r\,.
 \end{aligned}
\]
All in all, we conclude the upper energy-dissipation estimate 
\[
\calE_k(t,Z(t))+\calD(Z(0),Z(t)) \leq \calE_k(0,Z(0))+\int_0^t\partial_t\calE_k(r,Z(r))\,\mathrm{d}r \qquad \text{for all } t \in [0,T].
\]
 \par
 The lower estimate $\geq $
  follows from the stability condition \eqref{stab} via a well-established Riemann-sum technique, cf.\ e.g.\ \cite[Prop.\ 2.1.23]{MieRou-book} for a general result.  This concludes the proof that $Z$ is an Energetic solution to the  adhesive rate-independent system $(\Spz,\calE_k, \calD)$. 
\QED
\par 
Prior to proving Thm.\ \ref{thm:4.3}, we show in the following result that, for every fixed $t\in [0,T]$, 
the energies $\calE_k(t,\cdot)$ for the adhesive systems $\Gamma$-converge as $k\to\infty$ to the energy $\calE_\infty(t,\cdot)$ driving the brittle system 
\emph{with respect to the weak$^*$-topology of $ \SBV(\Omega;\{0,1\})$} (in the sense of \eqref{sense}).
Namely, we shall prove that 
\begin{subequations}
\label{Gamma-conv_E_k}
\begin{align}
&
\label{Ginf}
\text{$\Gamma$-$\liminf$ estimate:} && z_k \weaksto z \text{ in } \SBV(\Omega;\{0,1\}) \ \Rightarrow \ \liminf_{k\to\infty}\calE_k(t,z_k) \geq \calE_\infty(t,z),
\\
&
\label{Gsup}
\text{$\Gamma$-$\limsup$ estimate:} && \forall\, \hat{z} \in  \SBV(\Omega;\{0,1\}) \ \exists\, (\hat{z}_k)_k\subset 
 \SBV(\Omega;\{0,1\})\, : \ \hat{z}_k \weaksto \hat{z} \text{ in } \SBV(\Omega;\{0,1\})  \\
 \nonumber
 & && 
\qquad \qquad  \text{ and }  \limsup_{k\to\infty} \calE_k(t,\hat{z}_k) \leq \calE_\infty(t,\hat{z})\,.
\end{align}
\end{subequations}
\begin{lemma}
\label{PropGamma}
For every $t\in [0,T]$ the functionals $(\calE_k(t,\cdot)_k$ defined by 
 \emph{\eqref{ENADH}}
 $\Gamma$-converge as $k\to\infty$, in the sense of \emph{\eqref{Gamma-conv_E_k}}, to $\calE_\infty(t,\cdot)$ from 
  (\ref{ENBRI}).
\end{lemma}

\noindent {\it Proof.} To start with, observe that the upper estimate \eqref{Gsup} can be concluded by choosing the constant sequence $(\hat{z}_k: = \hat{z})_k$ 
as a recovery sequence. 
\par To verify the lower estimate \eqref{Ginf},  consider a sequence 
$z_k \weaksto z$ in  $\SBV(\Omega;\{0,1\})$ and 
 the corresponding sets 
of finite perimeter $ Z_k, Z$.  We distinguish two cases:
\begin{compactenum}
\item
If $\calL^d( Z\cap F(t))=0,$ then, by the positivity of 
 $\calJ_k(t, \widetilde Z_k)$ 
 and the lower semicontinuity of the
 perimeter w.r.t.\ 
strong $L^1$-convergence 
of  characteristic functions, 
 we find that 
 \[
 \liminf_{k\to\infty}\calE_k(t,z_k)
\geq P(Z,\Omega)=\calE_\infty(t,z).
\]
\item
Assume now that $\calL^d(Z\cap F(t))=c>0$, so that $\calJ_\infty(t,z)=\infty$.  
Since $z_k\to z$ strongly in $L^1(\Omega),$ we also have 
$f(t) z_k\to f(t) z$ strongly in $L^1(\Omega),$ and hence, for every $\eps>0$ we find an index 
$k_\eps\in\N$ such that for all $k\geq k_\eps$ we have that 
$\|f(t)z_k-f(t) z\|_{L^1(\Omega)}\leq\eps$. This implies that,   for every $\eps\in (0,c)$ there holds
\begin{equation*}
\liminf_{k\to\infty}\calJ_k(t,z_k)
=\liminf_{k\to\infty} \int_{\Omega} k f(t) z_k  \dd x 
\geq\liminf_{k\to\infty} k(c-\eps)=\infty,  
\end{equation*}
whence again \eqref{Ginf}.
\end{compactenum}
\QED
\par
We are now in a position to carry out the \underline{\textbf{proof of Theorem \ref{thm:4.3}:}}
Let us consider a sequence  $(Z_k)_k$  with $Z_k\in \SMM k$ for every $k\in \N$ and the associated sequence of characteristic functions $(z_k)_k$. 
Since the energy bound \eqref{energy-bound-interp} holds for a constant uniform w.r.t.\ $k\in \N$ as well,  and it is  inherited by the time-continuous limit, 
we infer
\[
\exists\, C>0 \ \forall\, k \in \N\, : \quad  \sup_{t\in [0,T]}\calE_k(t,z_k(t)) \leq C,  \ \text{ as well as } \  \mathrm{Var}_{\calD}(z_k; [0,T]) \leq C
\]
as the sequence $(z_k)_k$ is nonincreasing in time.
Hence we may repeat the very same compactness arguments as in the proof of Thm.\ \ref{thm:4.1} and conclude that there exist a (not relabeled) subsequence and $z\in L^\infty(0,T;\SBV(\Omega;\{0,1\})) \cap \BV([0,T];L^1(\Omega))$ such that \eqref{sense} holds, as well as 
\begin{equation}
\label{convk-L}
z_k\weaksto z \text{ in } L^\infty(0,T;\SBV(\Omega;\{0,1\})), \qquad
z_k\to z \text{ in } L^q ((0,T){\times} \Omega) 
\text{ for every } q \in [1,\infty).
\end{equation}
Thanks to Lemma \ref{PropGamma} we have 
\begin{equation}
\label{liminf_cons}
\liminf_{k\to\infty} \calE_k(t,z_k(t)) \geq \calE_\infty(t,z(t)), \quad \text{in particular } \calJ_\infty(t,z(t))=0, \quad \text{for all } t \in [0,T]. 
\end{equation}
\par
The limit passage as $k\to\infty$  in the stability condition \eqref{stab}, for $t\in (0,T]$ fixed,  again relies on the mutual recovery sequence condition, i.e.\ on the fact  that for every $\tilde z \in  \SBV(\Omega;\{0,1\})$ (with 
$  \calD( z(t),\tilde{z}) <\infty $ and $\calJ_\infty (t,\tilde z)=0$ to avoid trivial situations) there exists $(\tilde{z}_k)_k \subset  \SBV(\Omega;\{0,1\})$ such that 
\begin{equation}
\label{MRS-k}
\limsup_{k\to\infty} \left( \calE_k(t,\tilde{z}_k) {-}\calE_k(t,z_k(t)) {+} \calD(z_k(t),\tilde{z}_k) \right) \leq  \calE_\infty(t,\tilde{z}) {-}\calE_\infty(t,z(t)) {+} \calD( z(t),\tilde{z}) \,. 
\end{equation}
To this end, we resort to a construction completely analogous to that in \eqref{recovery-j}
and set
\begin{equation}
\label{recovery-k}
\tilde{z}_k: = \tilde z \chi_{A_k} +z_k(t) (1{-} \chi_{A_k})  \quad \text{with } A_k: = \{ x \in \Omega\, : \ 0\leq \tilde z(x) \leq z_k(t,x) \} \,.
\end{equation}
Exploiting 
the convergence properties \eqref{conv-dissip}, \eqref{limsup-perims}, as well as  the fact that 
\[
\limsup_{k\to\infty} \left( \calJ_k(t,\tilde{z}_k) -\calJ_k(t,z_k(t)) \right) = \limsup_{k\to\infty} \int_\Omega k f(t,x)\left( \tilde{z}_k(t,x) {-} z_k(t,x)\right) \ddx x 
\leq 0 = \calJ_\infty (t,\tilde z)-\calJ_\infty (t,z(t))
\]
(since $\tilde z_k \leq z_k$ a.e.\ in $\Omega$ by construction), we
obtain \eqref{MRS-k}. This proves that $Z\in \SMM \infty$. 
\par
In order to conclude that $Z$ is an Energetic solution to the brittle system $(\Spz,\calE_\infty,\calD)$ under the additional monotonicity assumption $\partial_t f \geq0$ and compatibilty condition \eqref{compatibility}, we first establish the upper energy-dissipation estimate
\begin{equation}
\label{enbal-infty-upper}
\calE_\infty(t,z(t))+\calD(z(0),z(t)) \leq \calE_\infty(0,Z_0) \quad \text{for all } t \in [0,T]
\end{equation}
by passing to the limit in  \eqref{enbal} as $k\to\infty$.
With this aim, we
observe that 
\begin{equation}
\label{convJk}
\calJ_k(t,z_k(t))=\int_\Omega kf(t)z_k(t)\,\mathrm{d}x\to0 \quad\text{as } k\to\infty \quad \text{for all } t \in [0,T].
\end{equation}
To check \eqref{convJk}, first of all  we apply the construction of the recovery sequence \eqref{recovery-k} 
to  $\tilde z:=z(t)$ to find a sequence $(\tilde z_k)_k$, associated with sets $(\widetilde Z_k)_k$.
As observed in the proof of Theorem \ref{thm:4.1} (cf.\ \eqref{limsup-perims}),  
the construction ensures 
in particular that $\limsup_{k\to\infty}\big(|\mathrm{D} \tilde{z}_k|(\Omega){-} |\mathrm{D} z_k(t)|(\Omega) \big)\leq |\mathrm{D} z(t)|(\Omega)- |\mathrm{D} z(t)|(\Omega) =0$.
Moreover, 
 $\tilde{z}_k \leq  \min\{ z_k(t), \, z(t) \}$ a.e.\ in $\Omega$ and 
$\tilde z_k\to z(t)$ in $L^1(\Omega)$, so that 
\[
\lim_{k\to\infty}\calD(z_k(t),\tilde z_k)=0.
\]
Also observe that $\calJ_k(t,\tilde{z}_k)=0$ for all $k\in\N,$ because $\calJ_\infty(t,z(t))=0$ 
by 
\eqref{liminf_cons}
 and because 
 $\widetilde Z_k\subset Z(t)$ by construction.  
Since $\calJ_k(t,z_k(t))\geq0$ for all $k\in\N$, 
these observations allow us to conclude from the stability of $z_k(t)$ that 
\begin{equation*}
0\leq \liminf_{k\to\infty} \calJ_k(t,z_k(t)) \leq\limsup_{k\to\infty}\calJ_k(t,z_k(t))
\leq\limsup_{k\to\infty}\big(P(\widetilde Z_k,\Omega)-P(Z_k(t),\Omega)\big) 
+\lim_{k\to\infty}\calD(z_k(t),\tilde Z_k)=0\,,
\end{equation*}
which gives \eqref{convJk}. 
\NEW We will now show that 
\begin{equation}
\label{power2zero}
\int_0^t \partial_t \calE_k(s,z_k(s)) \dd s = \int_{0}^{t}\int_{\Omega}k\partial_t f(s)z_{k}(s)\,\mathrm{d}x\,\mathrm{d}s \to 0  \quad\text{as } k\to\infty \quad \text{for all } t \in [0,T].
\end{equation}
Indeed,
using the assumption $\partial_t f \geq0$ the energy balance \eqref{enbal} for each $k\in\N$ gives
\begin{equation}
\label{BRAVA}
\begin{aligned}
0&\leq\int_{0}^{t}\int_{\Omega}k\partial_t f(s)z_{k}(s)\,\mathrm{d}x\,\mathrm{d}s
\\ & =P(Z_{k}(t),\Omega)+\calJ_{k}(t,z_{k}(t))+\calD(Z_{0}^{k},Z_{k}(t))-P(Z_{0}^{k},\Omega)  - \calJ_k(0,Z_0^k) 
\\
&\stackrel{(1)}{\leq} \calD(Z_{k}(t),\emptyset)
+\calD(Z_{0}^{k},Z_{k}(t))-P(Z_{0}^{k},\Omega)
\stackrel{(2)}{=}\calD(Z_{0}^{k},\emptyset)-P(Z_{0}^{k},\Omega)\,, 
\end{aligned}
\end{equation}
where {\footnotesize (1)} follows from observing that $ \calJ_k(0,Z_0^k) \geq 0$ and from choosing $\widetilde{Z} = \emptyset$  in the stability condition \eqref{stab} (which gives $P(Z_{k}(t),\Omega)+\calJ_{k}(t,z_{k}(t)) \leq \calD(Z_{k}(t),\emptyset)$), while {\footnotesize (2)} ensues from the fact that 
\[
\calD(Z_{k}(t),\emptyset) +\calD(Z_{0}^{k},Z_{k}(t))= \dissparam \calL^d(Z_k(t)) +\dissparam \calL^d(Z_0^k{\setminus}Z_k(t)) = \dissparam \calL^d(Z_0^k) = \calD(Z_{0}^{k},\emptyset)\,.
\]
We then observe that  
\[
\lim_{k\to\infty} \left( \calD(Z_{0}^{k},\emptyset){-}P(Z_{0}^{k},\Omega) \right)
\stackrel{(3)}{=} \calD(Z_{0},\emptyset)-P(Z_{0},\Omega) \stackrel{(4)}{=}0 \qquad \text{as } k \to \infty\,,
\]
where {\footnotesize (3)} is due to condition \eqref{conseq-well-prep} and {\footnotesize  (4)} follows from the compatibility condition \eqref{compatibility}.   
Hence we  conclude \eqref{power2zero}. Thus,
 \eqref{enbal-infty-upper} follows from passing to the limit as $k\to \infty$ in the upper estimate '$\leq$' of \eqref{enbal-infty}: the left-hand side is dealt with by lower semicontinuity arguments, while the limit passage on the right-hand side follows from the well-preparedness \eqref{well-prep} 
 of the initial data,
 and from \eqref{power2zero}. 
\par
The lower energy-dissipation estimate, 
i.e.\ the converse of inequality \eqref{enbal-infty-upper},
 follows from testing the stability condition at $t=0$ with $\widetilde Z = Z(t)$, which gives
\[
\calE_\infty(0,Z_0) \leq \calD(Z(0),Z(t)) + \calE_\infty(0,Z(t))  \leq    \calD(Z(0),Z(t)) + \calE_\infty(t,Z(t)).   
\]
Indeed, $\calE_\infty(0,Z(t)) - \calE_\infty(t,Z(t))= \calJ_\infty(0,Z(t)) - \calJ_\infty(t,Z(t))=0$,  since $Z(t) \cap F(t) =\emptyset$ and $F(0)\subset F(t) $  imply $Z(t) \cap F(0) =\emptyset$. 
\par
We have thus shown that $Z$ complies with the stability condition \eqref{stab} and with the energy-dissipation balance \eqref{enbal-infty}. 
This concludes the proof of Theorem \ref{thm:4.3}.
\QED
\begin{remark}
\label{rmk:compat4adh}
\upshape
Assume that the external loading $f$ for the \emph{adhesive} system complies with \eqref{power-control}, and that the initial datum $Z_0$ satisfies the compatibility  condition \eqref{compatibility}.
From \eqref{BRAVA} we then deduce that 
\[
0=\int_{0}^{t}\int_{\Omega}k\partial_t f(s)z_{k}(s)\,\mathrm{d}x\,\mathrm{d}s = \int_0^t \partial_t \calE_k(s,z_k(s)) \dd s
\]
for every $t\in (0,T]$ and for all $k\in \N$. 
Therefore, $\partial_t f z_k =0$ a.e.\ in $\Omega\times (0,T)$, i.e.\ the sets $Z_k$ fulfill
\begin{equation}
\label{weak-constraint}
Z_k(t) \subset ({\supp}(\partial_t f(t)))^c \qquad \foraa\, t \in (0,T)\,.
\end{equation}
Taking into account that $ (\supp(f(t)))^c \subset(\supp
(\partial_t f(t)))^c$,
\eqref{weak-constraint} may be interpreted as a weak form of the brittle constraint \eqref{brittle-constr-intro}. 
\end{remark}


\bibliographystyle{alpha}

\end{document}